\documentclass[3p,final]{elsarticle}

\usepackage{amssymb}
\usepackage{amsmath,mathtools}

\usepackage{dsfont}
\usepackage{hyperref}
\usepackage{cleveref}
\usepackage{amsthm}
\usepackage{numprint}
\usepackage{babel}[english]
\usepackage{caption}
\usepackage{subcaption}
\usepackage{ gensymb }
\usepackage{lineno}
\modulolinenumbers[5]
\usepackage{ragged2e}
\usepackage[dvipsnames]{xcolor}
\usepackage{microtype}
\usepackage{natbib}
\usepackage{graphicx}
\usepackage{multicol}
\usepackage{multirow}
\usepackage{relsize}
\usepackage{enumitem}
\usepackage{derivative}
\usepackage{booktabs}
\usepackage{bm}
\usepackage{soul}  
\DeclareMathOperator*{\argmin}{argmin} 
\DeclarePairedDelimiter{\abs}{\lvert}{\rvert}
\DeclarePairedDelimiter{\norm}{\lVert}{\rVert}

\newtheorem{definition}{Definition}

\begin{document}

\begin{frontmatter}



\title{Mapping locally nondominated curves in multiobjective topology optimization of compliance and volume}


\author[Diepenbeek,FlandersMake]{Tom De Weer}
\ead{tom.deweer at kuleuven.be}
\author[Diepenbeek,FlandersMake]{Elke Deckers}

\affiliation[Diepenbeek]{organization={KU Leuven Campus Diepenbeek, Department of Mechanical Engineering},
            addressline={Wetenschapspark 27}, 
            city={Diepenbeek},
            postcode={3590}, 
            country={Belgium}}

\affiliation[FlandersMake]{organization={Flanders Make@KU Leuven},
            country={Belgium}}

\begin{abstract}
This paper studies the nondominated set for topology optimization of compliance and volume.\
In the multiobjective topology optimization literature, scalarization is a popular method to approximate this set, which is better known as the Pareto frontier.\
However, previous work indicates that the pointwise approximation obtained by scalarization obfuscates its underlying mathematical structure: under an imposed length scale, it consists of locally nondominated curve segments.\
To investigate this structure, the paper uses a simple methodology to (i) generate local optima with a tightly controlled topology and (ii) use continuation to extend these points to locally nondominated curves.\
The paper then performs two rounds of numerical experiments.\
The first round applies the methodology to several cantilever examples.\
This leads to novel insights regarding topological complexity and symmetry and confirms prior conclusions regarding continuity, smoothness and convexity.\
The second round investigates two numerical instabilities: (i) a lack of uniqueness and (ii) discontinuities of the locally nondominated curves.\
The former is linked to an observed flatness of the optimization landscape, whereas the latter is attributed to branch splitting and linked to bifurcation theory.\
Namely, we classify these discontinuities as higher-dimensional analogues of the subcritical, symmetry-breaking pitchfork bifurcation.\
The paper concludes with suggestions for future research to address these features of the multiobjective optimization landscape.\
\end{abstract}




\end{frontmatter}

\section{Introduction}

Since the early work by Bendsøe and Kikuchi~\cite{bendsoe_kikuchi_1988}, topology optimization has matured into a powerful framework for generating optimized designs~\cite{SigMau2013,DeaGra2014}.\
Classically, the design question is formulated as an optimization problem with a single objective under one or more constraints.\
However, many design problems involve trade-offs between multiple objectives.\
In the topology optimization literature, this is often addressed with $\varepsilon$-constraint scalarization: the objective is minimized for various constraint bounds.\
For example, one of the oldest structural optimization problems considers minimization of compliance under a varying maximum volume constraint.\
Then, due to the nonconvexity of these scalarized subproblems, each one produces a local optimum.\
The result is a pointwise approximation that provides little insight into relationships between the local optima.\
The goal of this manuscript is to address this limitation by mapping the local optima that collectively define the solution structure of the multiobjective topology optimization problem.\\

Multiobjective optimization is concerned with finding the set of nondominated or Pareto-optimal points.\
Put simply, a nondominated point is one that cannot be improved for one objective without worsening another.\
When objectives are in conflict, many such points exist and together form the nondominated set, also known as the Pareto frontier.\
Thus, whereas single-objective optimization yields a single optimum, multiobjective optimization seeks the nondominated set.\
The final choice as to which nondominated solution is preferred is left to the user.\\

From the outset, the study of multiobjective topology optimization is complicated by the globality of the nondominance property.\ 
That is, a nondominated solution has no point in the \emph{entire} feasible set that improves upon it for all objectives.\ 
Given the large design freedom offered by topology optimization, verifying nondominance is generally intractable, akin to finding global optima in single-objective topology optimization.\
In practice, the topology optimization community relies on methods that produce local optima.\
As a result, \textit{the true shape of the nondominated set remains an open question in topology optimization}.\
Specifically, it remains unclear whether the Pareto frontier is continuous, smooth and convex.\
Recent research~\cite{DeWeer2026}, however, shows that the frontier \emph{can} be discontinuous, nonsmooth and nonconvex.\
It remains in question whether this possibility actually manifests.\\

To answer these questions, the Pareto frontier is commonly approximated with \emph{scalarization}.\
That is, the multiobjective problem is transformed (``scalarized'') into a series of single-objective problems, each of which is solved separately with a single-objective optimization routine that produces a local optimum.\
The two most popular methods are $\varepsilon$-constraint scalarization~\cite{deLeon2015,Ferrari2019,cool_TO_VA} and weighted-sum scalarization~\cite{ChenWu1998,Sato2017,Wang2022,Christensen2023,Chen2024,Almeida2025,Wotten2026}.\
Theoretically, the former can yield all nondominated solutions by a suitable variation of the scalarization parameters, whereas the latter is incapable of producing points in the nonconvex regions of the (upper image of the) nondominated set, provided standard assumptions from multiobjective optimization theory~\cite{Miettinen1998,Ehrgott2005,Eichfelder2008,DeWeer2026}.\\

\begin{figure}[h]
\centering
\begin{subfigure}[t]{0.5\linewidth}
  \centering
  \includegraphics[width=0.95\textwidth]{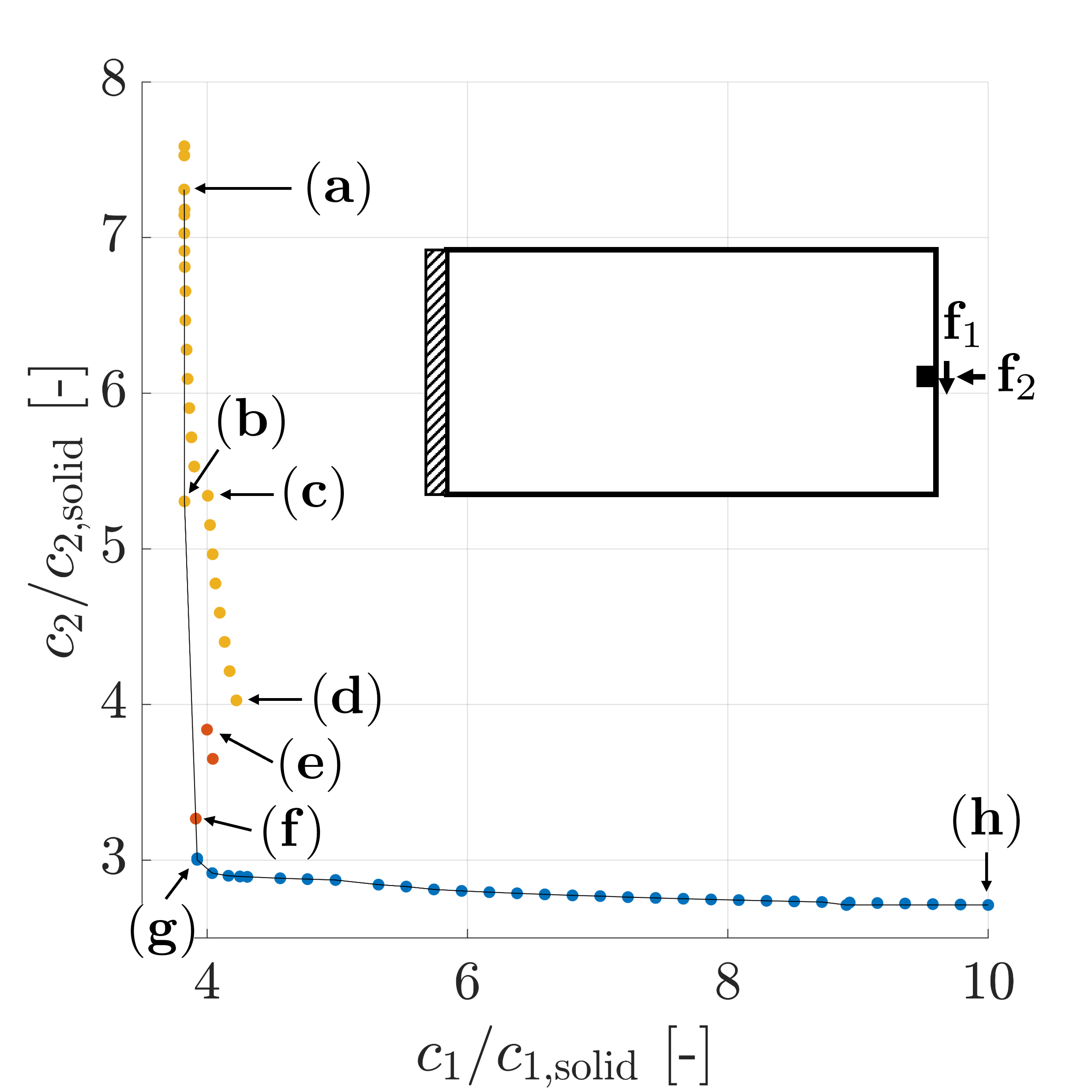}
\caption{Locally nondominated points.}
\end{subfigure}%
\begin{subfigure}[t]{0.5\linewidth}
  \centering
  \includegraphics[width=0.95\textwidth]{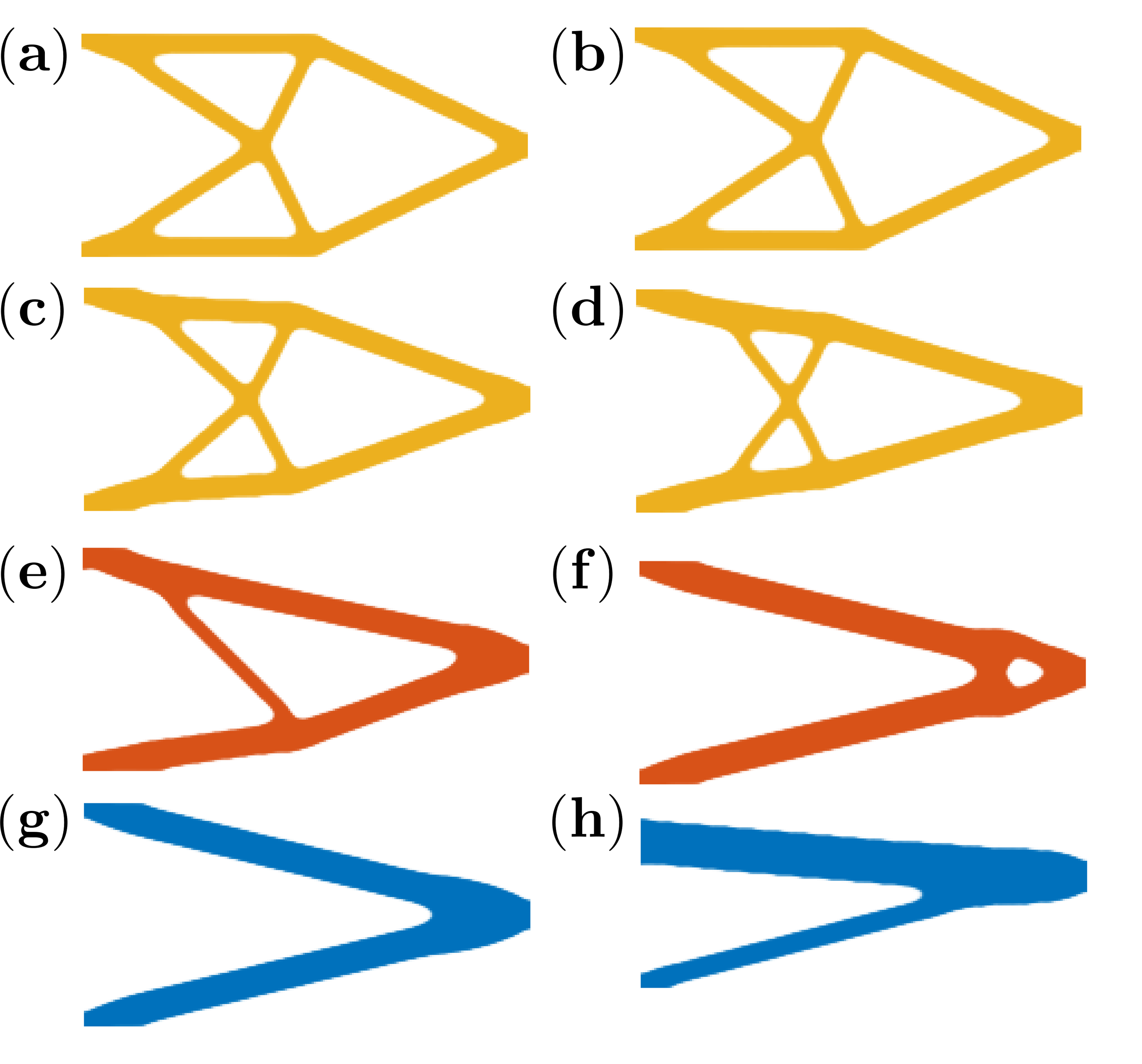}
\caption{Corresponding designs.}
\end{subfigure}%
\caption{Results of two $\varepsilon$-constraint scalarization variants applied to compliance minimization under two load cases (see top right of subfigure (a)).\
The blue, orange and yellow colors are used to indicate distinct topologies (1, 2 and 4 holes, resp.) between which a smooth design transformation is prohibited by the length scale.\
Results taken from \cite{DeWeer2026}, numerical details in \ref{app:TO}.\
The thin black line connects the nondominated points and approximates the nondominated set as a continuous curve.\
This is often done in the literature but can be misleading: the nondominated set has a gap between (b) and (f) and the corresponding designs have a different topology.\
}
\label{fig:scalarization}
\end{figure}

Scalarization yields a pointwise approximation: each local optimum is plotted as a point in the objective space.\
For an illustrative example, consider the approximation in \Cref{fig:scalarization} obtained by solving scalarized subproblems of a problem with two compliance minimization objectives under a maximum volume constraint.\
As is often done in the literature, the nondominated points are connected with a thin black line to indicate the extent of the nondominated set.\
However, this line is misleading for two reasons.\
First, one could assume that the nondominated set is a continuous curve.\
In fact, however, there is a large gap between the performances of designs (b) and (f), since designs (c-e) are dominated.\
There is no guarantee that nondominated points exist in this region of the objective space.\
Second, while continuity in the objective space might hold, the solid black line also implies that this continuity holds in the design space: it envisions a smooth transformation from design (b) to (f).\
Unfortunately, this is fundamentally impossible: design (b) has a different topology than design (f), but the imposed length scale prevents smooth topological changes.\
In other words, the two additional holes in design (b) cannot get arbitrarily small and hence must disappear discontinuously during any transformation from design (b) to (f).\\

A recent study~\cite{DeWeer2026} concludes that the nondominated set must be composed of segments of locally nondominated curves, each corresponding to a distinct topology.\
This mirrors the single-objective case, where the global optimum is a subset of the local optima.\
In \Cref{fig:scalarization}, at least three distinct locally nondominated curves can be observed: the blue, orange and yellow points correspond to 1-hole, 2-hole designs and  4-hole designs, respectively.\
\textit{The main goal of this paper is to better understand the observed piecewise nature of the nondominated set.}
To do so, the paper uses continuation and hole nucleation techniques to trace locally nondominated curves and explore the optimization landscape.\
To manage scope, the analysis is restricted to two objectives: compliance and volume minimization.\\

This work is related to several methods in the literature.\
These belong to a broader class of ``tracing methods'', as they aim to trace the frontier from one end to the other via continuation.\
Suresh~\cite{Suresh2010} uses topological sensitivity to continuously transform the design while moving along the frontier.\
During this procedure, the length scale is not enforced: the methodology produces designs with very small holes.\
Cesarano et al.~\cite{Cesarano2025} use a homotopy method for shape optimization and apply it to Pareto front tracing.\
Recently, Gangl and Winkler~\cite{Gangl2026} applied the homotopy approach to topology optimization.\
Their methodology starts from a local optimum of a relaxed multiobjective problem and extends this point with homotopy to obtain a series of stationary points.\
The combination of (i) relaxation to allow for grey design regions and (ii) the search for stationary points instead of minimizers allows transitions between local minima with different topology.\ 

This work is distinct from the previous tracing methods in two respects.\
First, it aims to tightly control the length scale.\
Penalization and filtering are used to prevent large grey regions containing fictitious material.\
Furthermore, any observed changes in the topology are used to terminate the branch.\
Second, it seeks to answer fundamental questions regarding the shape of the true frontier and therefore focuses on repeated continuation.\
That is, by mapping multiple branches it aims to approximate the global frontier instead of only producing designs with a single topology.\
From this point of view, \emph{the proposed methodology is a tracing method that gradually improves the approximation of the nondominated set by successively generating topologically distinct branches.}\\

After providing a general theoretical background in \Cref{sec:theoretical_background}, \Cref{sec:methodology,sec:numerical_examples_1,sec:numerical_examples_2} each contain one of the three main contributions of this work:
\begin{itemize} 
\item First, \Cref{sec:methodology} provides a straightforward methodology to map locally nondominated curves, based on (i) continuation of an $\varepsilon$-constraint scalarization and (ii) a hole nucleation heuristic.\
Due to its simplicity, it is an ideal exploratory tool for both the present study and for researchers in other subfields of topology optimization.
\item Second, \Cref{sec:numerical_examples_1} applies the methodology to several cantilever examples to investigate (locally) nondominated curves.\
In particular, \Cref{sec:numerical_examples_1} studies (i) the convexity and continuity of locally nondominated curves, (ii) how they combine to form the nondominated set and (iii) the (a)symmetry and complexity of the corresponding designs.\
\item Third, \Cref{sec:numerical_examples_2} investigates and explains two numerical instabilities observed in \Cref{sec:numerical_examples_1}:\
\begin{itemize}
\item[a)] The non-uniqueness of the obtained curves is linked to a local flatness of the optimization landscape.
\item[b)] Discontinuities in the locally nondominated curves are linked to subcritical, symmetry-breaking pitchfork bifurcations.
\end{itemize}
Adressing these instabilities, which stem from fundamental properties of the optimization landscape, is a necessity for further development of tracing methods.\
\end{itemize}
\Cref{sec:conclusions} then summarizes the conclusions of this paper.\

\section{Theoretical background and definitions} \label{sec:theoretical_background}

The purpose of this section is to provide a theoretical background, establish the relevant terminology and define the concepts required to clarify the goal and motivation of this work.\
\Cref{sec:MOO_background} introduces optimality in multiobjective optimization by providing definitions for (local) efficiency and nondominance.\
Then, \Cref{sec:scalarization_background,sec:continuation_background} present two primary approaches for solving multiobjective optimization problems: scalarization and continuation.\
Next, \Cref{sec:connectivity_background} formalizes the central goal of this paper by defining the locally efficient and nondominated curves in the context of topology optimization.\
Finally, \Cref{sec:bifurcation_background} establishes the connection between continuation-based methods and bifurcation theory.\
It also introduces the pitchfork bifurcation, which serves as a low-dimensional analogue to some of the the behaviour observed in this work.\

\subsection{Multiobjective optimization} \label{sec:MOO_background}

The below section provides a brief introduction to key concepts in multiobjective optimization.\
For an extended overview, tailored to the topology optimization community, the reader is referred to the overview in \cite{DeWeer2026}.\
More information can be found in reference works such as \cite{Miettinen1998,Ehrgott2005,Eichfelder2008}.\\

\begin{figure}
     \centering
      \includegraphics[width=0.8\textwidth]{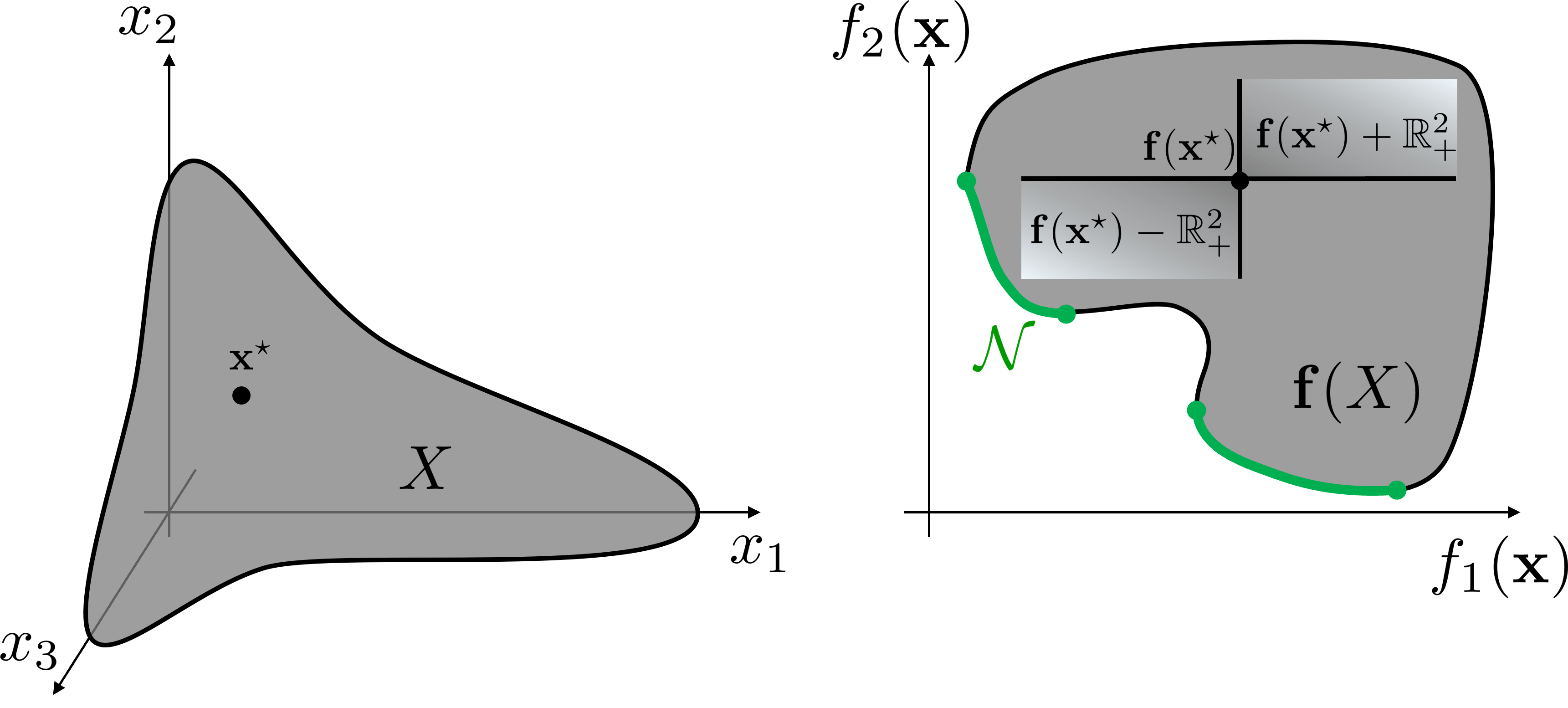}
     \caption{Optimization variable space and feasible set (left) and image space (right). The points in green show the nondominated set. The illustrated solution $\mathbf{x^\star}$ is not efficient since $\mathbf{f}(\mathbf{x}^\star)$ has feasible points in its bottom left quarter plane.}
    \label{fig:MOO_illustration}
 \end{figure}

Consider a multiobjective optimization problem with $m$ objectives
\begin{equation}
  \min_{\mathbf{x} \in X} \, \mathbf{f}(\mathbf{x}) = (f_1(\mathbf{x}),\dots,f_m(\mathbf{x}))^\top,
  \label{eq:mop}
  \tag{MOP($\mathbf{f},X$)}
\end{equation}
where $X \subseteq \mathbb{R}^n$ is the $n$-dimensional feasible set and $f_i : \mathbb{R}^n \to \mathbb{R}$ are the objective
functions.\
This work considers $m=2$ objectives (volume and compliance) but \ref{eq:mop} is stated for the general multiobjective case.\
\Cref{fig:MOO_illustration} illustrates the optimization variable space and the image space.\\

A solution $\mathbf{x}^\star$ is \emph{efficient} when it cannot be improved for one objective without worsening another.\
Its image $\mathbf{f}(\mathbf{x}^\star)$ is then called \emph{nondominated}.\
Mathematically, this is captured by the following two definitions:
\begin{definition} \label{def:efficient}
A solution $\mathbf{x}^\star \in X$ is \emph{efficient} if there exists no $\mathbf{x} \in X$ with $f_i(\mathbf{x}) \leq f_i(\mathbf{x}^\star)$ for all $i \in \{ 1, \ldots, m\}$ and $f_j(\mathbf{x}) < f_j(\mathbf{x}^\star)$ for at least one $j \in \{ 1,\ldots, m\}$.
\end{definition}
\begin{definition} \label{def:nondominated}
A point $\mathbf{f}(\mathbf{x}^\star) \in \mathbf{f}(X)$ is nondominated if $\left( \mathbf{f}(\mathbf{x}^\star) - \mathbb{R}_{+}^m \right) \cap \mathbf{f}(X) = \mathbf{f}(\mathbf{x}^\star)$.
\end{definition}
Here, $\mathbb{R}_{+}^m:= \left\{ \mathbf{y} \in \mathbb{R}^m: y_i \geq 0 \right\}$ is the nonnegative orthant of $\mathbb{R}^m$.\
\Cref{fig:MOO_illustration} illustrates $\mathbf{f}(\mathbf{x}^\star) \pm \mathbb{R}^m$ for a dominated point $\mathbf{f}(\mathbf{x}^\star)$.\\

All efficient points are collected into the \emph{efficient set} $\mathcal{E}(\mathbf{f}, X)$ of \ref{eq:mop} and its image $\mathcal{N}(\mathbf{f},X)=\mathbf{f}(\mathcal{E})$ is called the \emph{nondominated set}.\
For conciseness, these are often shortened to $\mathcal{E}$ and $\mathcal{N}$.\
Note that in engineering contexts, the nondominated set is often called the \emph{Pareto front(ier)}.\
For clarity, this work will only use the term nondominated set.\\

Obtaining the efficient and nondominated sets, or at least a ``good'' approximation of them, is the overarching goal of multiobjective optimization.\
However, topology optimization is notorious for its many local optima.\
Global optimality is only achievable and provable in rare, simplified cases~\cite{Dalklint2024,Zhou2025}.\
A straightforward application of \Cref{def:efficient,def:nondominated} is thus complicated, as efficiency and nondominance rely on global properties (cf. ``no $\mathbf{x} \in X$'' in \Cref{def:efficient}).\
To remain preciese, local efficiency and nondominance are now defined below.\

\begin{definition} \label{def:locally_efficient}
A solution $\mathbf{x}^\star \in X$ is \emph{locally efficient} if there exists an $\epsilon > 0$ such that there is no $\mathbf{x} \in X \bigcap B_n(\mathbf{x}^\star,\varepsilon) $ with $f_i(\mathbf{x}) \leq f_i(\mathbf{x}^\star)$ for all $i \in \{ 1, \ldots, m\}$ and $f_j(\mathbf{x}) < f_j(\mathbf{x}^\star)$ for at least one $j \in \{ 1,\ldots, m\}$.
\end{definition}
\begin{definition} \label{def:locally_nondominated}
A point $\mathbf{f}(\mathbf{x}^\star) \in \mathbf{f}(X)$ is locally nondominated if $\mathbf{x}^\star$ is locally efficient.
\end{definition}
The main difference is between global and local efficiency/nondominance  is that the latter's definition only holds in a small ball (radius $\epsilon$) around $\mathbf{x}^\star$.\
Finally, denote $\mathcal{E}_{\mathrm{loc}}(\mathbf{f},\mathbf{X})$ to be the set of all locally efficient solutions, or $\mathcal{E}_{\mathrm{loc}}$ for short.\
Similarly, $\mathcal{N}_{\mathrm{loc}}(\mathbf{f},\mathbf{X}) = \mathbf{f}(\mathcal{E}_{\mathrm{loc}})$, or $\mathcal{N}_{\mathrm{loc}}$ for short.\\

\begin{figure}
     \centering
      \includegraphics[width=0.8\textwidth]{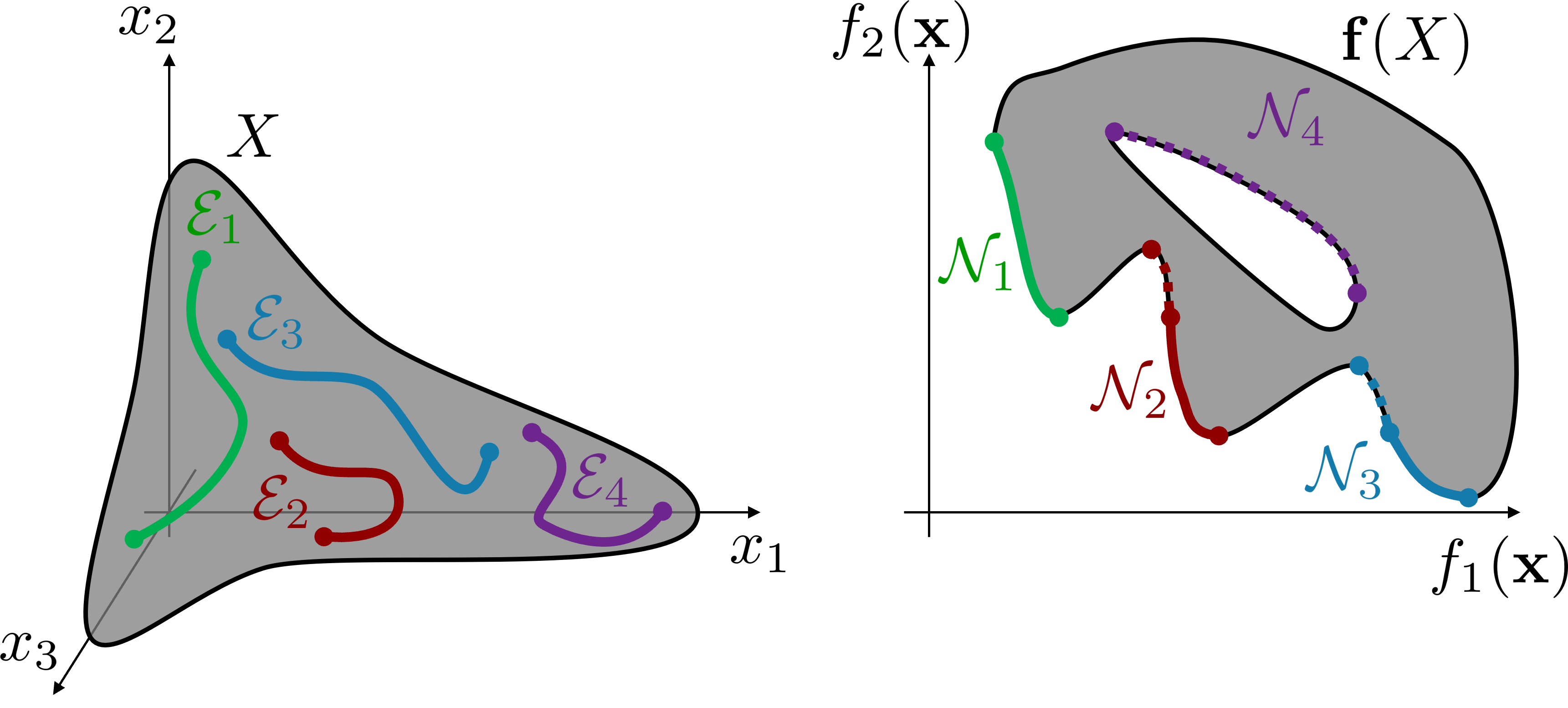}
     \caption{Illustration of four locally efficient and nondominated curves.\
     The dashed portions of the locally nondominated curves reflect points which are not globally nondominated.\
     Thus, each point in $\mathcal{N}_4$ is completely dominated.\
     Instead, the globally nondominated set consists of pieces of $\mathcal{N}_1$, $\mathcal{N}_2$ and $\mathcal{N}_3$.\
     }
    \label{fig:LMOO_illustration}
 \end{figure}

\Cref{fig:LMOO_illustration} illustrates the concept of locally efficient and nondominated curves.\
The main idea, further elaborated for in \Cref{sec:connectivity_background}, is that a continuity of the design variables and objective functions leads to locally efficient solutions lying on curves.\
These locally efficient curves then map to locally nondominated curves under $\mathbf{f}$.\

\subsection{Scalarization} \label{sec:scalarization_background}

Scalarization is an established technique to solve multiobjective optimization problems.\
The idea is to transform \ref{eq:mop} into a parameter-dependent family of single-objective problems

\begin{equation}
  \min_{\mathbf{x} \in X^\star (\boldsymbol{\mu})} \, f^\star(\mathbf{x},(\boldsymbol{\mu})),
  \label{eq:mop_mu}
 \tag{$\mathrm{SOP}(\boldsymbol{\mu})$}
\end{equation}

where the new feasible set $X^\star \subseteq \mathbb{R}^n$ and the scalarized objective function $f^\star \in \mathbb{R}^n \times \mathbb{R}^k \rightarrow \mathbb{R}$ depend on a parameter vector $\boldsymbol{\mu} \in \mathbb{R}^p$, with $p$ the number of parameters.\
This work refers to the parameters $\boldsymbol{\mu}$ as the \emph{scalarization parameters}.\\

By way of example, the well-known weighted-sum scalarization is defined as follows.\
\begin{definition}
The weighted-sum scalarization is an instance of \ref{eq:mop_mu} where $\boldsymbol{\mu}=\mathbf{w} \in \mathbb{R}^m_+$ (so $p=m$) , $X^\star=X$ and $f^\star (\mathbf{x}) = \mathbf{w}^\top \mathbf{f}(\mathbf{x})$, with $\mathbf{w}$ a weight vector such that $\sum_{i=1}^m w_i = 1$.\
It can be written as 
\begin{equation} \label{eq:WS_formulation}
  \min_{\mathbf{x} \in X} \, \mathbf{w}^\top \mathbf{f}(\mathbf{x}). 
  \tag{$\mathrm{WS}(\mathbf{w})$}
\end{equation}
\end{definition}

Another well-known example is the $\varepsilon$-constraint scalarization, defined as follows.\
\begin{definition}
The $\varepsilon$-constraint scalarization is an instance of \ref{eq:mop_mu} where $\boldsymbol{\mu}=\left( k, \boldsymbol{\varepsilon}^\top \right)^\top$, with $k \in \left\{ 1, \ldots, m \right\}$ an index to select one of the objectives and $\boldsymbol{\varepsilon} \in \mathbb{R}^{m-1}$ containing constraint values on the other objectives.\
More succinctly, $f^\star = f_k$ and $X^\star = X \, \cap \, \left\{ \mathbf{x}: f_i(\mathbf{x}) \leq \varepsilon_i, \, \forall i \in I \right\}$ where $I=\left\{ 1, \ldots, m \right\} \setminus k$.\
It can be written as
\begin{equation} \label{eq:EC_formulation}
  \min_{\mathbf{x} \in X} f_k(\mathbf{x}) 
  \quad \text{s.t.} \quad f_i(\mathbf{x}) \leq \varepsilon_i \quad \forall i \in \left\{ 1, \ldots, m \right\} \setminus k.
    \tag{$\mathrm{EC}(\bm{\varepsilon})$}	
\end{equation}
\end{definition}
The $\varepsilon$-constraint scalarization is used in this work to find locally nondominated points.\\

An important property of a scalarization technique is the recovery property: can the nondominated set $\mathcal{N}$ be obtained (``recovered'') through a suitable variation of the scalarization parameters?\
In the case of \ref{eq:WS_formulation}, $\mathcal{N}$ is recoverable if $\mathbf{f}(X)+\mathbb{R}_+^m$ is a convex set~\cite{Ehrgott2005}.\
Convexity is a sufficient but quite strict condition.\
Recent research~\cite{DeWeer2026} indicates that $\mathbf{f}(X)+\mathbb{R}_+^m$ is not convex in topology optimization, namely in regions where locally nondominated curves intersect.\
This is confirmed by the numerical results of this work.\\

Practically, if \ref{eq:WS_formulation} is applied to a nonconvex \ref{eq:mop}, i.e., to an \ref{eq:mop} where $\mathbf{f}(X)+\mathbb{R}_+^m$ is a nonconvex set, then the concave part is not recoverable since these are not global minimizers of \Cref{eq:WS_formulation}.\
In contrast, $\mathcal{N}$ is recoverable with \ref{eq:EC_formulation} through a suitable variation of $k$ and $\bm{\varepsilon}$.\
Furthermore, note that a weighted-sum objective can still be employed to approximate concave portions of the frontiers when either (a) \Cref{eq:WS_formulation} is not solved to global but only local optimality, as in this work, or (b) one is not restricted to local minimizers but allows for stationary points.\
The latter is used by e.g. Hillermeier~\cite{Hillermeier2001} and discussed in the next section.\

\subsection{Continuation} \label{sec:continuation_background}

The concept of finding locally nondominated points through continuation goes back to the original work by Rakowska et al. in 1991~\cite{Rakowska1991} and was followed by several decades of improvements~\cite{Lundberg1993,Hillermeier2001,Guddat2007,Harada2007,Potschka2011,Ringkamp2012}.\
It has been applied to chemical engineering~\cite{Seferlis1996} and, more recently, large-scale dynamical systems~\cite{Keler2017}.\
The core idea is to formulate the multiobjective conditions of optimality as a set of parameterized nonlinear equations 
\begin{equation} \label{eq:gxmu}
\mathbf{g}(\mathbf{x}, \bm{\lambda})=\mathbf{0},
\end{equation}
where $\bm{\lambda}$ is a vector collecting one or more unspecified parameters.\
Starting from a solution $\mathbf{x}^\star$ that solves \Cref{eq:gxmu} for $\bm{\lambda}^\star$, a continuation then finds subsequent points by using $\mathbf{x}^\star$ as the initial guess for finding a solution of \Cref{eq:gxmu} with $\bm{\lambda}^\star + \Delta \bm{\lambda}$.\
As continuation methods resemble the solution of an ordinary differential equation (ODE), they can be solved with predictor-corrector methods~\cite{Potschka2011}.\
The activation of constraints can also cause bifurcations~\cite{Martin2016}.\\

Instead of directly exploiting multiobjective conditions of optimality, such as the Fritz-John conditions in \cite{Martin2016}, one can also use scalarization.\
That is, locally nondominated points can be found by varying the scalarization parameters in \ref{eq:mop_mu} (so $\bm{\lambda}=\bm{\mu}$).\
Note that this idea can extend the applicability of the weighted-sum scalarization: by continuing local optima with respect to $\mathbf{w}$, concave regions of the nondominated set can be explored.\
A practical consequence, however, is that classic optimizers are no longer suitable: instead of finding minimizers of \ref{eq:mop_mu}, it requires stationary points of \ref{eq:mop_mu}, with stationary expressed as $\mathbf{g}(\mathbf{x}, \bm{\lambda})=\mathbf{0}$.\
Consequently, non-minimizing stationary points of \Cref{eq:WS_formulation}, such as saddle points and local maxima, must also be considered.\
For example, the recent work of Gangl and Winkler~\cite{Gangl2026} uses a Newton method to find and continue stationary points.\
To avoid these complications, this work instead applies continuation to an $\varepsilon$-constraint scalarization.\
This is further detailed in \Cref{sec:continuation}.\\

\subsection{Connectivity and topology} \label{sec:connectivity_background}

One of the core assumptions of a continuation method is that two subsequent solutions $\mathbf{x}_1$ and $\mathbf{x}_2$, which solve \Cref{eq:gxmu} for $\bm{\mu}_1$ and $\bm{\mu}_2$, are connected.\
That is, there exists a continuous map $\gamma: [0,1] \rightarrow \mathbb{R}^n$ such that
\begin{itemize}
\item  $\gamma(0)=\mathbf{x}_1,\gamma(1)=\mathbf{x}_2$; and
\item $\gamma(t)$ is locally efficient for all $t \in [0,1]$.
\end{itemize}
This idea is now extended to topology optimization, where the optimization variables $\mathbf{x}$ encode a twodimensional shape $\Omega \in \mathbb{R}^2$.\
In addition to requiring $\gamma(t)$ to be locally efficient, solutions in topology optimization are only connected if the corresponding shapes are homeomorphic.\
In practice, this requires that they have the same number of holes.\
This is now formalized.\
\begin{definition} \label{def:local_frontier}
Consider a subset of locally efficient solutions $\mathcal{E}^\star \subseteq \mathcal{E}_{\mathrm{loc}}$, where every element $\mathbf{x}$ of $\mathcal{E}^\star$ represents a shape $\Omega$.\
If for all pairs $\mathbf{x}_1 \in \mathcal{E}^\star$ (corresponding to $\Omega_1$) and  $\mathbf{x}_2\in \mathcal{E}^\star$ (corresponding to $\Omega_2$), there exists a  continuous map $\gamma: [0,1] \rightarrow X$ such that
\begin{itemize}
\item  $\gamma(0)=\mathbf{x}_1,\gamma(1)=\mathbf{x}_2$; and
\item $\gamma(t) \in \mathcal{E}^\star$ for all $t \in [0,1]$; and
\item  for all $t_a, t_b \in [0,1]$, $\Omega_a$ corresponding to $\mathbf{x}_a=\gamma(t_a)$ is homeomorphic to $\Omega_b$ corresponding to $\mathbf{x}_b=\gamma(t_b)$;
\end{itemize}
then $\mathcal{E}^\star$ is a locally efficient curve and its image $\mathcal{N}^\star=\mathbf{f}(\mathcal{E}^\star)$ is a locally nondominated curve.\
\end{definition}

A major insight of recent research~\cite{DeWeer2026} is that the efficient set in multiobjective topology optimization 
\begin{itemize}
\item consists of  designs $\Omega$ with different topologies; and
\item contains locally efficient curves, as defined in \Cref{def:local_frontier}.\
\end{itemize}
However, the conclusions in~\cite{DeWeer2026} are obfuscated due to scalarization.\
The approximations are pointwise and, since the optimization process cannot control the topology, each point can lie on a distinct local frontier.\
Nonetheless, locally efficient and nondominated curves are observed and \emph{the methodological goal of this work is to map locally efficient and nondominated curves in topology optimization by maintaining strict control of topology and length scale.}\

\paragraph{Terminology} For conciseness, a pair $\mathcal{E}^\star,\mathcal{N}^\star$ of a locally efficient curve and its corresponding nondominated curve is from now on denoted as a \emph{branch}.\
Multiple branches together form a \emph{tree}, and the goal of this work is to study trees by mapping their branches.\

\subsection{Bifurcation theory} \label{sec:bifurcation_background}

The numerical experiments in this work lead to several branches that exhibit inflection points and discontinuities.\
That is, whereas most locally nondominated curves are continuous and convex, a small minority are not.\
Further investigation, carried out in \Cref{sec:bifurcations}, shows that most of these inflections are due to bifurcations.\
Local minima transform into local maxima, i.e., they go from stable to unstable, and the optimizer moves away from the unstable optimum to a ``nearby'', stable local minimum.\
As long as this local minimum corresponds to a design with the same topology, the branch mapping continues and leads to nonconvex locally nondominated curves.\\

The bifurcations observed in this work are related to symmetry-breaking subcritical pitchfork bifurcations.\
This term is now briefly explained.\
For a more detailed overview of bifurcation theory, the reader is referred to reference works such as the book by Strogatz~\cite{Strogatz2018}.\\

Bifurcation theory is classically the study of parameterized nonlinear dynamical systems

\begin{equation} \label{eq:dynamical_system}
\dot{\mathbf{x}}=\mathbf{g}(\mathbf{x}, \bm{\lambda}),
\end{equation}
 where $\mathbf{x}$ and $\bm{\lambda}$ are the state and parameter vectors and $\mathbf{g}$ is a generic function mapping $(\mathbf{x}, \bm{\lambda})$ to the state derivative vector $\dot{\mathbf{x}}$.\
Note that $\dot{\mathbf{x}}$ is the \emph{time} derivative of $\mathbf{x}$.\
A bifurcation occurs when a property of the dynamical system changes under small deviations of $\bm{\lambda}$.\
For the purpose of this work, the relevant property is the stability of an equilibrium, i.e., a pair $(\mathbf{x}, \bm{\lambda})$  for which $\dot{\mathbf{x}}=\mathbf{g}(\mathbf{x}, \bm{\lambda})=\mathbf{0}$.\\

To stress the relation to optimization, we consider here an unconstrained, parameter-dependent single-objective optimization problem, which could be considered an instance of the multiobjective continuation problem studied in this work.\
Let $x \in \mathbb{R}$ be the optimization variable, $\lambda$ the parameter and $f: (\mathbb{R},\mathbb{R})  \mapsto \mathbb{R}$ the parameter-dependent objective function.\
That is, we consider the optimization problem
\begin{equation}
\min_{x} f(x,\lambda)
\end{equation}
A point $\mathbf{x}^{\star}$ is a stationary point when $\partial f(x^{\star},\lambda) / \partial x = 0$.\
It is a local minimum if $\partial^2 f(x^{\star},\lambda) / \partial x^2 > 0$ and a local maximum if $\partial^2 f(x^{\star},\lambda) / \partial x^2 < 0$.\
The equivalence to the dynamical system of \Cref{eq:dynamical_system} is straightforward: the path $x(t)$ of a ball that rolls down the $f(x,\lambda)$ landscape along the steepest gradient is described by the dynamical system $\dot{x} = - \frac{\partial f(x, \lambda)}{ \partial x}$.\
Thus, a dynamical system $\dot{x}=g(x,\lambda)$ is equivalent to minimization of $f(x, \lambda)=-\int g(x,\lambda) \,\mathrm{d}x$, irrespective of the integration constant.\\

\begin{figure}
     \centering
      \includegraphics[width=0.99\textwidth]{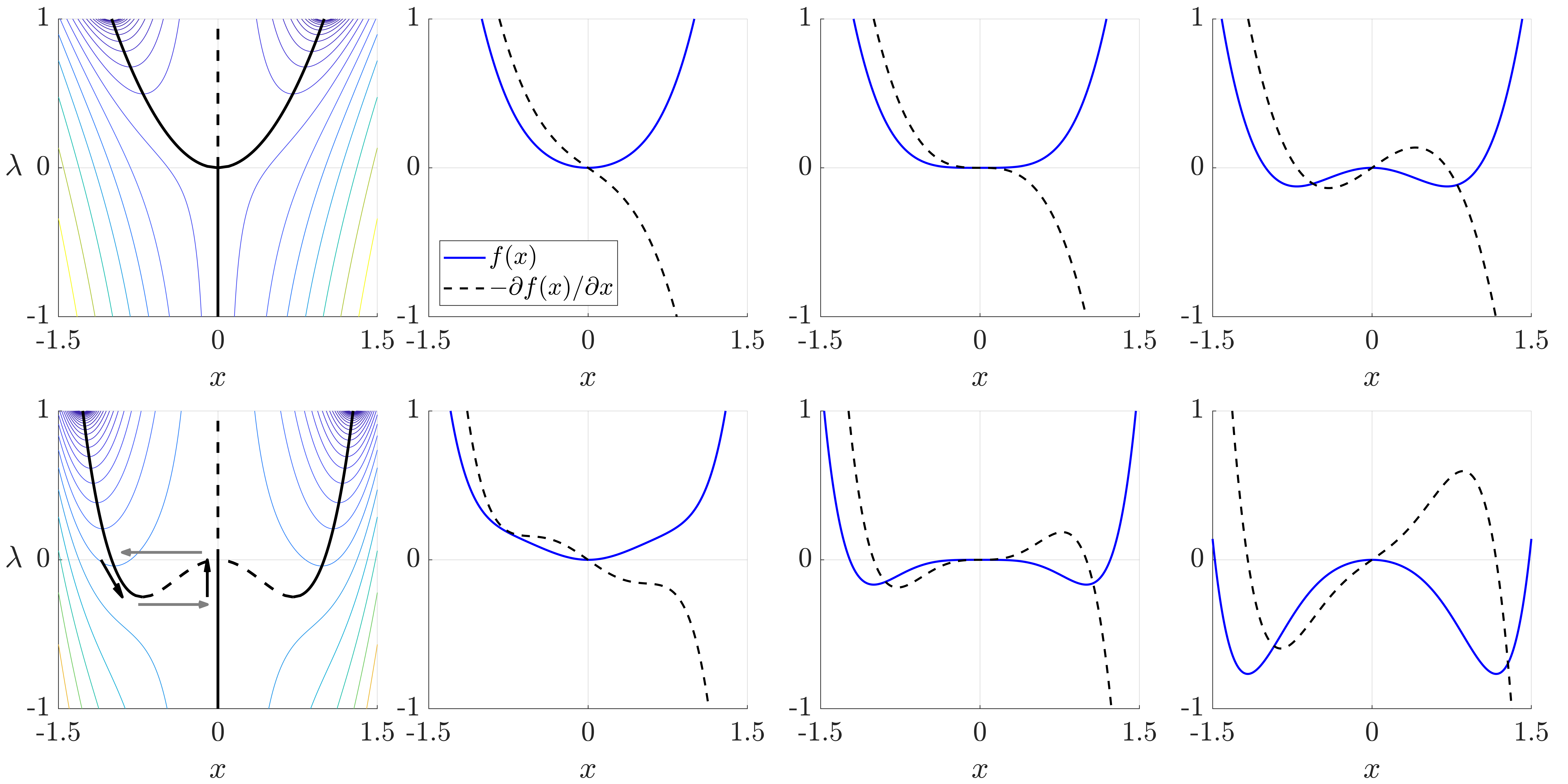}
     \caption{Illustration of two types of pitchfork bifurcation.\
      Top row: supercritical pitchfork.\
       Bottom row: subcritical pitchfork.\
       Left column: contour lines of $f(x, \lambda)$ in the $(x, \lambda)$ plane, with solid black lines denoting local minima and dashed lines denoting local maxima.\
       Right three columns: $f(x, \lambda)$ as function of $x$ for three values of $\lambda$: $\lambda=-1$ (left),  $\lambda=0$ (middle),  and $\lambda=1$ (right).\
       The bottom left figure, showing the $(x, \lambda)$ plane of the subcritical pitchfork bifurcation, has several arrows to indicate a hysteresis when varying $\lambda$.\ 
       Additionally, the grey arrows denote a sudden jump caused by a loss of stability\
     }
    \label{fig:pitchforks}
 \end{figure}

\Cref{fig:pitchforks} illustrates the two types of pitchfork bifurcations.\
The prototypical \emph{supercritical} pitchfork bifurcation occurs for $g(x, \lambda)=\lambda x - x^3 = \partial f(x^{\star},\lambda) / \partial x$.\
This corresponds to minimization of $f(x,\lambda) = -\lambda x^2 + \frac{x^4}{2}$, where the integration constant is omitted and a rescaling occurred.\
For $\lambda < 0$, $\min_x f(x,\lambda)$ has one minimum at $x=0$.\
For $\lambda > 0$, it has two local minima (at $x=\pm \sqrt{\lambda}$) and one local maximum at $x=0$.\
The change at $\lambda=0$, which transforms a minimum into a maximum and creates two new, symmetric minima, is known as a pitchfork bifurcation because of its visualization in the $(x,\lambda)$ plane.\\

The $-x^3$ term for the supercritical pitchfork is stabilizing, causing the local minimum at $x=0$ to transform smoothly into two new local minima.\
If the sign is reversed (-$x^3 \mapsto +x^3$), this stabilization mechanism is lost.\
Then, a small deviation will cause solutions of $\dot{x}=g(x,\lambda)$ to converge to a ``nearby'' equilibrium.\
A higher order stabilization term $x^5$ is added to model these nearby equilibria, resulting in the prototypical \emph{subcritical} pitchfork bifurcation, $\partial g(x,\lambda) / \partial x = \lambda x + x^3 - x^5$ or $f(x,\lambda) = -\lambda x^2 - \frac{x^4}{2} + \frac{x^6}{3}$, visualized in the bottom row of \Cref{fig:pitchforks}.\\

In the subcritical case, as $\lambda$ passes through zero, the stability of the $x=0$ minimum is lost and no new minimum is nearby.\
This leads to a sudden jump in $x$, which moves either left or right.\
Then, decreasing $\lambda$ again leads to another jump back to the minimum at $x=0$.\
Continuation on objective function landscapes with subcritical bifurcations can thus lead to jumps and hysteresis.\
Conversely, jumps and hysteresis during continuation can be used to identify subcritical bifurcations.\
This observation is used later in this work to identify subcritical bifurcations.\\

Both sub- and supercritical pitchfork bifurcations exhibit a symmetry in the optimization landscape: $f(x)$ is symmetric around the stationary point at $x=0$.\
After the bifurcation, the symmetric equilibrium loses stability and the optimum shifts to an asymmetric solution.\
Consequently, pitchfork bifurcations are commonly referred to as symmetry-breaking phenomena.\
In this work, they correspond to the transition of a symmetric to an asymmetric optimum.\

\section{Methodology} \label{sec:methodology}

To map branches in topology optimization, this work combines three methods, which are now introduced:
\begin{itemize}
\item \Cref{sec:TO}: A classical density-based topology optimization framework, employed at high $\beta$ projection to prevent spontaneous hole nucleation.\
\item \Cref{sec:continuation}: A continuation method, based on $\varepsilon$-constraint scalarization, to expand locally nondominated points to locally nondominated curves.\
\item \Cref{sec:hole_nucleation}: A hole nucleation heuristic, based on topological derivatives,  to explore branch candidates.\
\end{itemize}

\begin{figure}
     \centering
     \captionsetup{width=\textwidth}
      \includegraphics[width=0.7\textwidth]{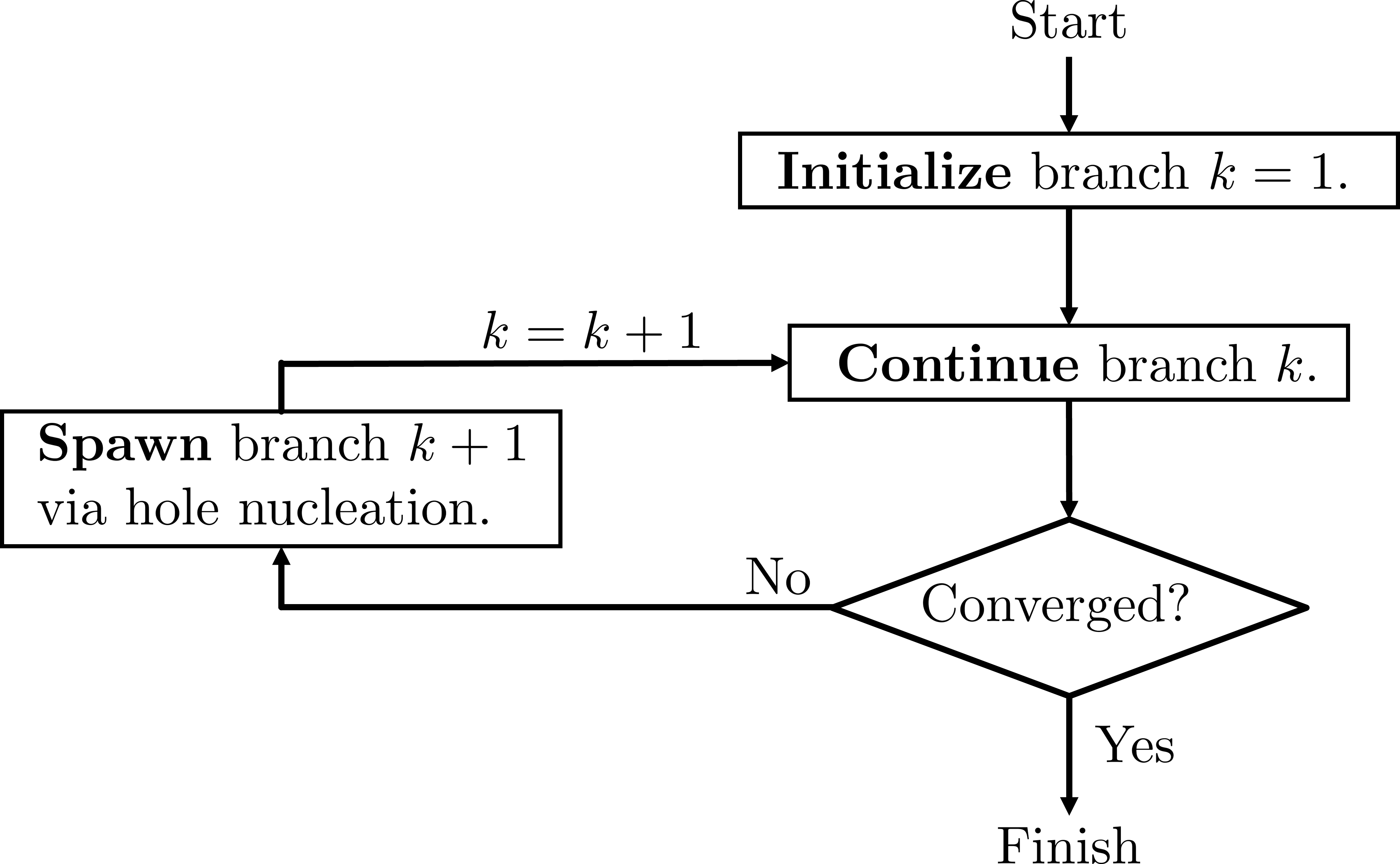}
     \caption{Overview of the continuation and hole nucleation routine.}
    \label{fig:workflow}
 \end{figure}

\Cref{fig:workflow} shows a schematic overview of the methodology.\
Given a starting design for the first branch, the method expands the current branch through continuation.\
Then, the method either converges or uses a hole nucleation heuristic to generate an initial guess for the next branch.\
This work employs a maximum number of iterations as convergence criterion but can converge sooner if no suitable new branch can be spawned.\

 \subsection{Topology optimization} \label{sec:TO}

This work considers topology optimization of compliance and volume.\
The main goal of the topology optimization scheme is to bring an initial guess to optimality (i.e., local nondominance), which is now explained.\
Further numerical details are in \ref{app:TO}.\

\subsubsection*{Design representation}

\begin{figure}
     \centering
     \captionsetup{width=0.6\textwidth}
      \includegraphics[width=0.95\textwidth]{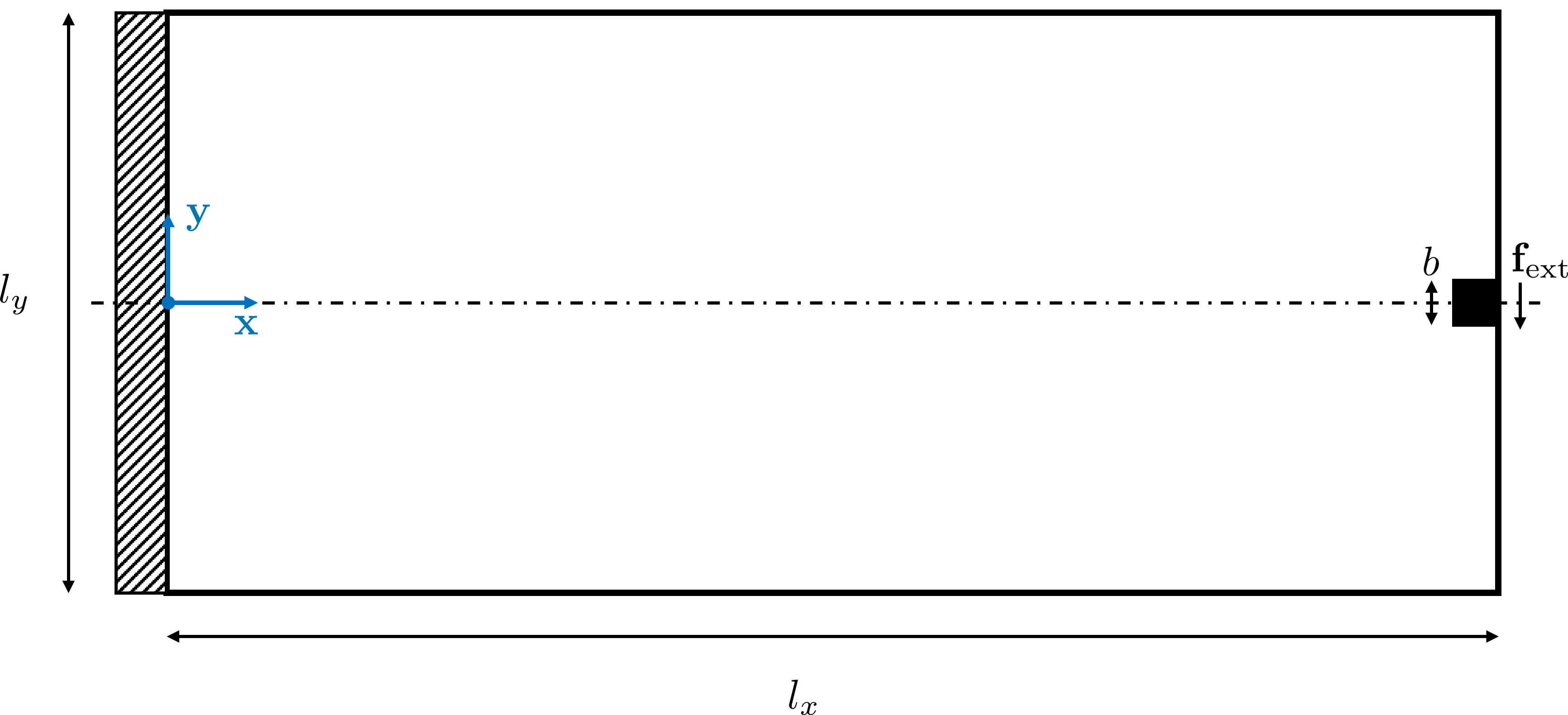}
     \caption{Cantilever case studied in this work.}
    \label{fig:Cantilever_setup}
 \end{figure}

The study considers a cantilever of size $l_x$ by $l_y$, shown in \Cref{fig:Cantilever_setup}.\
The topology optimization routine is density-based, and hence the cantilever is meshed with $n=n_x \times n_y$ square first-order finite elements.\
Then, the design is parameterized by the design vector $\mathbf{x} \in \mathbb{R}^n$, with each entry corresponding to a finite element such that $x_i=0$ denotes a void and $x_i=1$ denotes solid material.\\

The routine employs a three-field filtering scheme to prevent checkerboarding and to impose a length scale.\
It consists of the PDE filter from Lazarov and Sigmund~\cite{Lazarov2011}, with filter radius $R$,  followed by a smoothed Heaviside projection~\cite{Guest2004} with steepness $\beta$ and threshold $\eta$.\
The thresholding is performed at three $\eta$ levels, yielding three distinct designs: an eroded design ($\mathbf{x}_e$ at $\eta_e=0.6$), a blueprint design  ($\mathbf{x}_b$ at $\eta_e=0.5$) and a dilated design  ($\mathbf{x}_d$ at $\eta_d=0.4$).\
To avoid boundary effects, the consistent boundary conditions from Wallin et al.~\cite{Wallin2020} are employed.\
The three design instances (eroded/blueprint/dilated) are combined via the robust approach of Wang et al.~\cite{wang2011projection}.\
This enforces a length scale of $b_{\min}$ for both void and solid regions as $\beta \rightarrow \infty$~\cite{qian2013topological}.\\

The density-based, three-field filtering scheme is employed with a high $\beta$ projection value, specified in \ref{app:TO}, to ensure it acts like a level-set method.\
Namely, at high $\beta$, the design is almost black and white except for a thin layer of elements near solid-void boundaries.\
Because of the filtering, the design cannot represent shapes smaller than an induced length scale.\
Furthermore, because sensitivities are almost zero in solid and void regions, the spontaneous appearance of holes is strongly inhibited.\
This ensures tight control over the creation of holes.\

\subsubsection*{Compliance and volume}

The volume $V(\mathbf{x})$ of a design is computed as $V(\mathbf{x}) =  \sum_{i=1}^{n} x_{d,i} V_i$, with $V_i$ the volume of element $i$ and $x_{d,i}$ the value of the dilated density vector $\mathbf{x}_d$ corresponding to element $i$.\\

The compliance is computed as $C(\mathbf{x}) = \mathbf{u}^\top \mathbf{f}_{\mathrm{ext}}$, where $\mathbf{f}_{\mathrm{ext}}$ is the external force vector and $\mathbf{u}$ is the displacement vector.\
It is found by solving the state equations for linear elasticity:
\begin{equation}
\mathbf{K}(\mathbf{x}_e) \mathbf{u}=\mathbf{f}_{\mathrm{ext}}.
\end{equation}
Here, $\mathbf{K}(\mathbf{x}_e)$ is the stiffness matrix for the eroded design vector.\
The routine uses SIMP interpolation~\cite{Bendsoe1989} with a penalization constant equal to $3$ to discourage the appearance of grey material.\\

Both compliance and volume are normalized by their solid design values.\
That is, the compliance is defined as $c(\mathbf{x})=C(\mathbf{x}) / c_{\mathrm{solid}}$ and the volume as $v(\mathbf{x})=V(\mathbf{x}) / v_{\mathrm{solid}}$, with $c_{\mathrm{solid}}$ and $v_{\mathrm{solid}}$ the compliance and volume of the completely filled design.\
Furthermore, only the eroded compliance and dilated volume are computed since these are the most stringent.\
For conciseness, the compliance and volume functions are simply denoted as $v(\mathbf{x})$ and $c(\mathbf{x})$ and $\mathbf{x}$ is implicitly understood to refer to the correct, filtered design vector.\
When designs are plotted, the blueprint realization is used.\

\subsubsection*{Formulation}

This work considers the following optimization problem.\

\begin{equation}
\begin{aligned}
\min_{\mathbf{x} \in \left[ \mathbf{0}, \mathbf{1} \right]} \quad & \left( c(\mathbf{x}), \, v(\mathbf{x}) \right)^\top \\
\text{s.t.} \quad & c(\mathbf{x}) < c_{\max}
\end{aligned}
\label{eq:min_cv}
\end{equation}

The upper bound $c_{\max}$ prevents infinitely compliant structures.\
 
\subsection{Continuation} \label{sec:continuation}

Locally nondominated points are expanded to locally nondominated curves via a combination of scalarization and continuation.\

\subsubsection*{Scalarization}

The $\varepsilon$-constraint scalarization is used to find locally nondominated points of \Cref{eq:min_cv}.\
In practice, there are two variants: a compliance minimization variant

\begin{equation} \label{eq:minC_variant}
  \min_{\mathbf{x}} c(\mathbf{x}) 
  \quad \text{s.t.} \quad v(\mathbf{x}) \leq \varepsilon_v
    \tag{$\mathrm{EC_c}(\varepsilon_v)$}	
\end{equation}

and a volume minimization variant

\begin{equation} \label{eq:minV_variant}
  \min_{\mathbf{x}} v(\mathbf{x}) 
  \quad \text{s.t.} \quad c(\mathbf{x}) \leq \varepsilon_c \leq c_{\max} .
    \tag{$\mathrm{EC_v}(\varepsilon_c)$}	
\end{equation}

As stated, $\mathcal{N}$ can be recovered via a suitable selection of the variant and $\varepsilon$.\
This also holds in the local case: every locally nondominated point is recoverable, i.e., it represents a local optimum of one or both of the above variants.\
In this work, the recoverability property is exploited to find locally nondominated curves by solving successive instances of \ref{eq:minC_variant} and \ref{eq:minV_variant}, which is now explained further.\

\subsubsection*{Continuation of $k$ and $\bm{\varepsilon}$}

Consider an initial state $(\mathbf{x}_0, \boldsymbol{\mu}_0)$, with $\mathbf{x}_0$ corresponding to a shape $\Omega_0$, that solves \ref{eq:mop_mu}.\
In other words, $\mathbf{f}(\mathbf{x}_0)$ is locally nondominated.\
The goal of the continuation methodology is to find the corresponding branch, i.e., the locally efficient curve that goes through $\mathbf{x}_0$ and its corresponding locally nondominated curve.\
In practice, this is achieved by solving \ref{eq:mop_mu} for a suitable variation of $\boldsymbol{\mu}$ and using previous solutions as initial guesses for subsequent solves.\
The goal is then to obtain a set of locally efficient designs $\mathbf{x}_1$, $\mathbf{x}_2, \ldots$ that represent shapes $\Omega_1, \Omega_2, ...$ that are homeomorphic to $\Omega_0$ and each other.\\

Note that there is a twofold assumption underlying the methodology, namely that (i) there exists a locally efficient curve through $\mathbf{x}_0$ and that (ii) this curve is unique.\
The numerical experiments will provide no counterexample to the first assumption but do show a lack of uniqueness.\
This is further discussed in \Cref{sec:numerical_examples_2}, but already mentioned here since it constitutes one of the main conclusions of this paper: locally efficient curves can split.\\

\begin{figure}
     \centering
     \captionsetup{width=\textwidth}
      \includegraphics[width=0.5\textwidth]{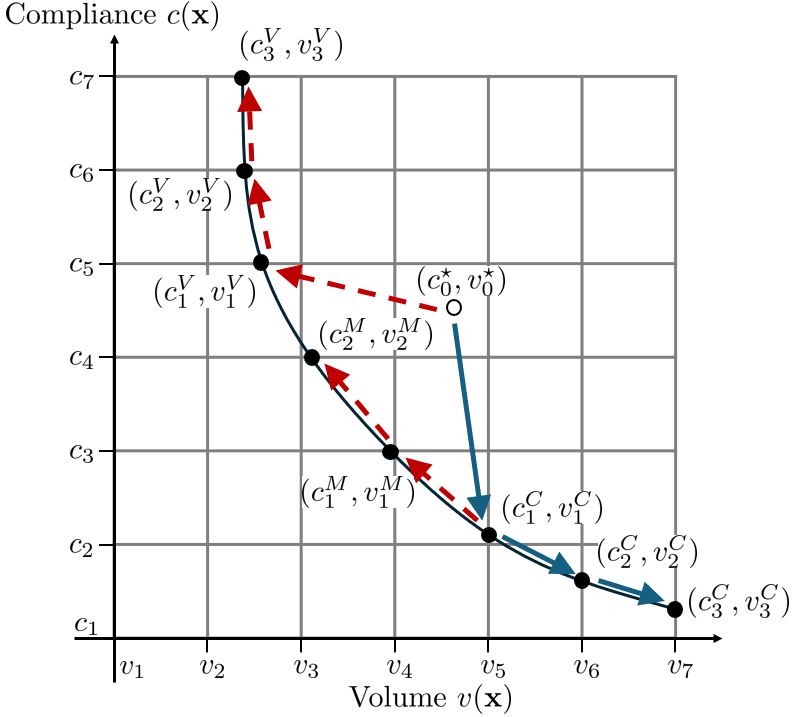}
     \caption{Illustration of the continuation workflow in the $(c,v)$ space. Filled circles denote locally nondominated points. Arrows denote volume and compliance minimization runs (dashed red arrow and full blue, respectively).}
    \label{fig:continuation_workflow}
 \end{figure}

\Cref{fig:continuation_workflow} illustrates the continuation workflow.\
The  $\left( c(\mathbf{x}), v(\mathbf{x}) \right)$ subspace $\left[ 1, c_{\max} \right] \times \left[ 0, 1 \right]$ is first discretized into a grid by dividing $\left[1, c_{\max} \right]$ into $n_c$ equispaced compliance values $c_1, c_2, \ldots, c_{n_c}$ and $\left[0, 1 \right]$ into $n_v$ equispaced volume values $v_1, v_2, \ldots, v_{n_v}$.\\

The continuation scheme starts from an initial guess $\mathbf{x}_0^\star$, with corresponding objective function values $c_0^\star$ and $v_0^\star$, generated with the hole nucleation heuristic discussed in \Cref{sec:hole_nucleation}.\
The continuation scheme consists of three subcontinuation runs.\\

The first subcontinuation brings $\mathbf{x}_0^\star$ to optimality and then moves it towards minimal compliance.\
This is done by solving \ref{eq:minC_variant} with $\mathbf{x}_0^\star$ as initial guess and $\varepsilon_v=v_i$, with $i=\argmin_{i \in \left\{ 1, \ldots, n_v\right \}} || v_i-v_0^\star || \, \mathrm{s.t.} \, v_i \geq v_0^\star$.\
That is, $i$ is such that $v_i$ is the smallest value still larger than $v_0^\star$.\
If no such $v_i$ exists, the search in the minimal compliance direction stops.\
Otherwise, the resulting solution is called $\mathbf{x}_1^C$ and corresponds to a locally nondominated point $(c_1^C, v_1^C)$.\
Then the procedure continues by using $\mathbf{x}_1^C$ as initial guess for a compliance minimization with $\varepsilon_v=v_{i+1}$ until $i+1>n_v$.\\

The second subcontinuation goes towards minimal volume.\
Starting with $\mathbf{x}_0^\star$ as initial guess, \ref{eq:minV_variant} is solved for increasing values of $\varepsilon_c$, starting with $\varepsilon_c=c_i$, with $i=\argmin_{i \in \left\{ 1, \ldots, n_c\right \}} || c_i-c_0^\star || \, \mathrm{s.t.} \, c_i \geq c_0^\star$.\
This yields locally efficient solutions $\mathbf{x}_1^V$, $\mathbf{x}_2^V, \ldots$ corresponding to locally nondominated points $(c_1^V, v_1^V), (c_2^V, v_2^V), \ldots$ until $i+1>n_c$.\\ 

If $x_0^\star$ is strongly nondominated, then $(c_1^C, v_1^C)$ and $(c_1^V, v_1^V)$ can be far apart.\
To close this gap, a third subcontinuation is done.\
Namely, successive volume minimization runs are initiated, starting from $\mathbf{x}_1^V$ as initial guess with $\varepsilon_c=c_i$ where $i=\argmin_i || c_i-c_0^\star || \mathrm{s.t.} c_i \geq c_1^V$.\
The resulting solution, denoted $\mathbf{x}_1^M$, is used as initial guess for increasing values of $\varepsilon_c$ until $c_{i+1}=c_1^V$.\\

Subcontinuations can end prematurely for two reasons.\
First, if a topological change is detected between design $\mathbf{x}_i^{C/V/M}$ and $\mathbf{x}_{i+1}^{C/V/M}$, then the latter is discarded and the subcontinuation ends.\
This often occurs for the first subcontinuation, towards minimal compliance, since if $c(\mathbf{x}) \rightarrow 1$ then $v(\mathbf{x}) \rightarrow 1$ and holes get filled.\
Second, if a solution $\mathbf{x}_{i}^{C/V/M}$ is ``close'' to an already found solution, it is discarded and the subcontinuation again ends.\
The detection of topological changes and closeness between two solutions is detailed in \ref{app:continuation}.\\ 

\subsection{Hole nucleation} \label{sec:hole_nucleation}

The first branch uses the fully solid design as initial guess.\
Subsequent initial guesses are found via a hole nucleation heuristic based on topological derivatives.\
The main idea behind the heuristic is to place holes with a diameter at the minimum length scale $b_{\min}$ in locations with low strain energy.\
The nucleation itself is done by setting the density value to $0$ for all elements in the ``hole region'', i.e., elements whose centroid lies within the circular region around a ``puncture location''.\\

Suitable puncture locations are those where the topological derivative for the compliance is low.\
The topological derivative encodes the effect (on the compliance) of puncturing the design at location $\mathbf{p} \in \mathbb{R}^2$ with an infinitesimally small hole.\
For the numerical examples in this work, the material's Poisson ratio is $\nu=1/3$ and thus the topological derivative for the compliance is proportional to the strain energy density.\
Therefore, this work uses a hole nucleation approach similar to the one used in~\cite{Christiansen2014,Andreasen2020} but extended to multiobjective optimization.\\

\begin{figure}[h]
     \centering
     \captionsetup{width=\textwidth}
      \includegraphics[width=0.5\textwidth]{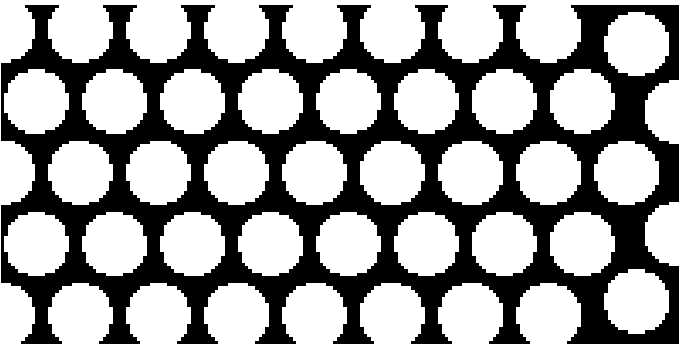}
     \caption{Hole candidates for a symmetric $2\,\times\,1$ cantilever setup. Filter radius is $15$ element widths, for a length scale of $0.095$.}
    \label{fig:hole_grid}
 \end{figure}

To restrict the number of hole nucleation points, only a preselected number of hole candidates are considered.\
These are generated in a way that depends on the design domain, length scale and boundary conditions.\
\Cref{fig:hole_grid} provides an example for the symmetric $2\,\times\,1$ cantilever studied in \Cref{sec:numerical_examples_1}.\

The hole nucleation scheme chooses a design and punctures it with at most $n_H$ holes.\
This is done as follows.\

The set of designs that can be punctured is first limited to those whose performance is currently nondominated.\
For each such design $k$, the hole nucleation metric
\begin{equation}
\phi_{i,k}= (\mathbf{u}_{i,k}^e)^\top \mathbf{K}_{i,k}^e \mathbf{u}_{i,k}^e,
\end{equation}
is computed for every hole candidate in the grid (e.g., for every hole candidate in \Cref{fig:hole_grid}).\
The index $i$ corresponds to the element at the center of a hole candidate.\
Thus, $\mathbf{u}_{i,k}^e$ and $\mathbf{K}_{i,k}^e$ are the elemental displacement vector and stiffness matrix of element $i$ for design $k$, respectively.\
Next, the hole nucleation metric is dedimensionalised.
\begin{equation}
\tau_{i,k} = \frac{\phi_{i,k} - \phi_{\min, k}}{\phi_{\max,k} - \phi_{\min,k}}
\end{equation}
Here, $\phi_{\min,k}$ and $\phi_{\max,k}$ are the minimum and maximum values of $\phi_{i,k}$ in design $k$.\
Due to the normalization, $0 \leq \tau_{i,k} \leq 1$.\\

\begin{figure}
     \centering
     \captionsetup{width=\textwidth}
      \includegraphics[width=\textwidth]{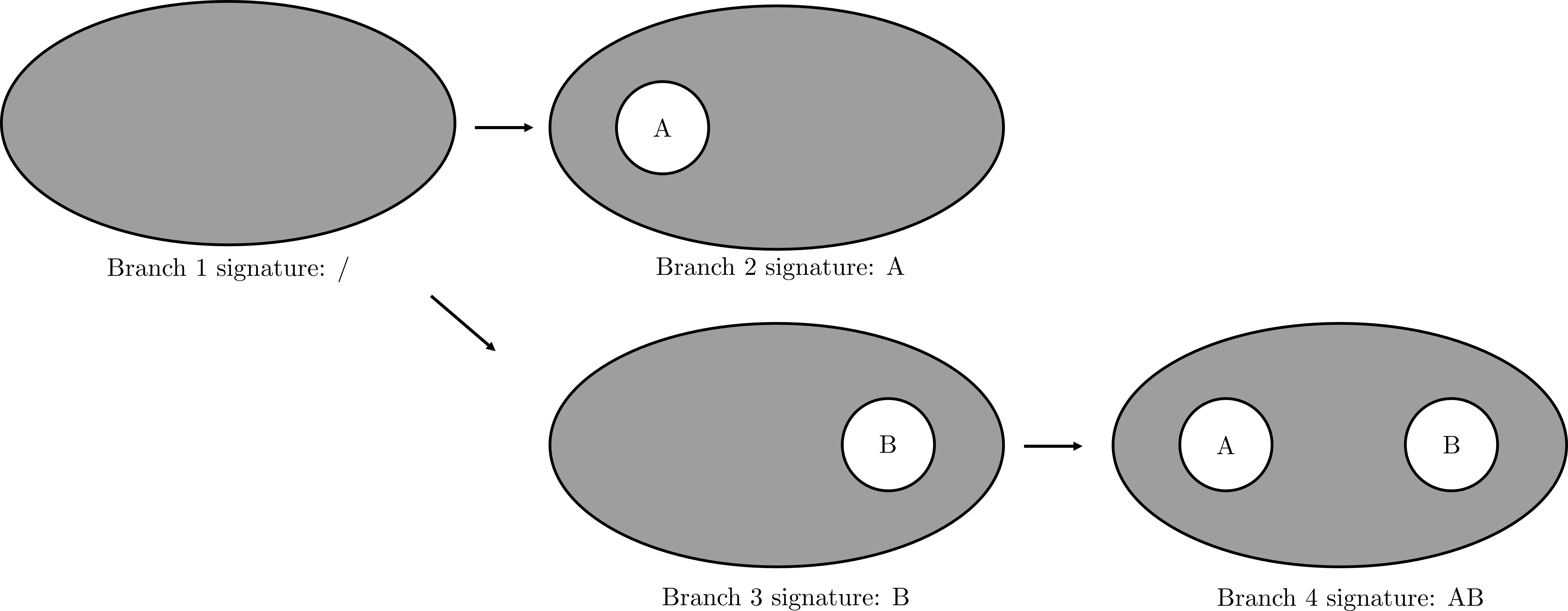}
     \caption{Illustration of a hole signature. The starting branch (branch 1), has no holes and hence an empty hole signature. Branch 2 is found by puncturing a design from branch 1 at hole location A. Thus, branch 2 has hole signature A. Similarly, branch 3 has hole signature B. Branch 4 is found by puncturing branch 3 at hole location A, yielding signature AB. The hole nucleation scheme prohibits the creation of a ``new'' branch 5 by puncturing branch 2 at B, since the resulting signature AB would coincide with branch 4.}
    \label{fig:puncture_map}
 \end{figure}

The hole nucleation metric provides a way to gauge which hole combination is likely to improve performance.\
However, before comparison three restrictions are applied.
\begin{itemize}
\item First,  holes must increase the genus of the structure.\
It is therefore checked if they do by applying the puncture and computing the number of holes with a flood fill algorithm (details in \ref{app:hole_nucleation}).\
\item Second, only holes with $\tau_{i,k} \leq \mathrm{tol}_{\tau}$ are considered.\
This prohibits the nucleation of holes in regions with a (relatively) high strain energy.\
\ref{app:TO} lists the value of $\mathrm{tol}_{\tau}$ for each numerical example.\
\item Third, puncturing a design with (at most $n_H$) holes is not allowed if this is likely to lead to an already existing design.\
This is done by comparing the hole signature (see \Cref{fig:puncture_map}).
\end{itemize}

After the above constraints are considered, the combination of holes for which the average $\tau_{i,k}$ is the lowest is selected.\
If there are no hole combinations that satisfy the above restrictions, the mapping procedure terminates.\

The proposed hole nucleation scheme is similar to the original bubble method by Eschenauer et al.~\cite{Eschenauer1994}.\
There, for a fixed maximum volume, compliance minimization is done with a shape optimization routine.\
Whenever convergence is reached, a hole (``bubble'') is added in a region with a low topological derivative.\
By repeating this procedure, the bubble method generates a series of design with the same volume but increasing genus and decreasing compliance.\
The method proposed in this work can be viewed as a multiobjective extension of the bubble method.\

\section{Numerical examples: mapping locally nondominated curves} \label{sec:numerical_examples_1}

This section applies the branch mapping methodology to several setups to study the nondominated set and its underlying locally nondominated curves.\
It contains four numerical studies, each of which studies the continuity and convexity of both the locally nondominated curves and their composition into the nondominated set.\
Naturally, since global optimality is not guaranteed, only an approximation of the nondominated can be discussed.\
\begin{itemize}
\item \Cref{sec:cantilever_beam_example} considers on a $2\,\times\,1$ cantilever beam loaded at the bottom right.\
It is meant primarily as an illustration of the mapping methodology.\
\item \Cref{sec:highC} considers the high compliance regime for the same cantilever, loaded at the bottom and middle right.\
It investigates and confirms that the lowest-topology design, a single beam, overtakes the two-bar design at sufficiently high compliance.\
\item \Cref{sec:symmetry_study_2x1} considers a cantilever beam of size $2\,\times\,1$, loaded at the middle right, and studies the symmetry of the solutions under the symmetric boundary conditions.\
\item \Cref{sec:symmetry_study_4x1} repeats the study of \Cref{sec:symmetry_study_2x1} for a $4\,\times\,1$ cantilever to investigate the effect of the length.\
\end{itemize}
A detailed overview of optimization parameters (mesh discretization, convergence tolerances, ...) can be found in the appendices.\
Namely, \ref{app:TO}, \ref{app:hole_nucleation},  \ref{app:continuation} and \ref{app:values} detail the topology optimization, hole nucleation, continuation and numerical parameters, respectively.\\

To aid the interpretation, many of the obtained branches are manually filtered.\
For example, branches with similar-looking designs are often removed to improve clarity.\
Their appearance signals a lack of uniqueness, which is further studied in \Cref{sec:numerical_examples_2}.\
A full accounting of omitted branches can be found in \ref{app:values}.\

\subsection{Cantilever beam} \label{sec:cantilever_beam_example}

\begin{figure}[h]
\centering
\begin{subfigure}[t]{0.33\linewidth}
  \centering
  \includegraphics[width=0.95\textwidth]{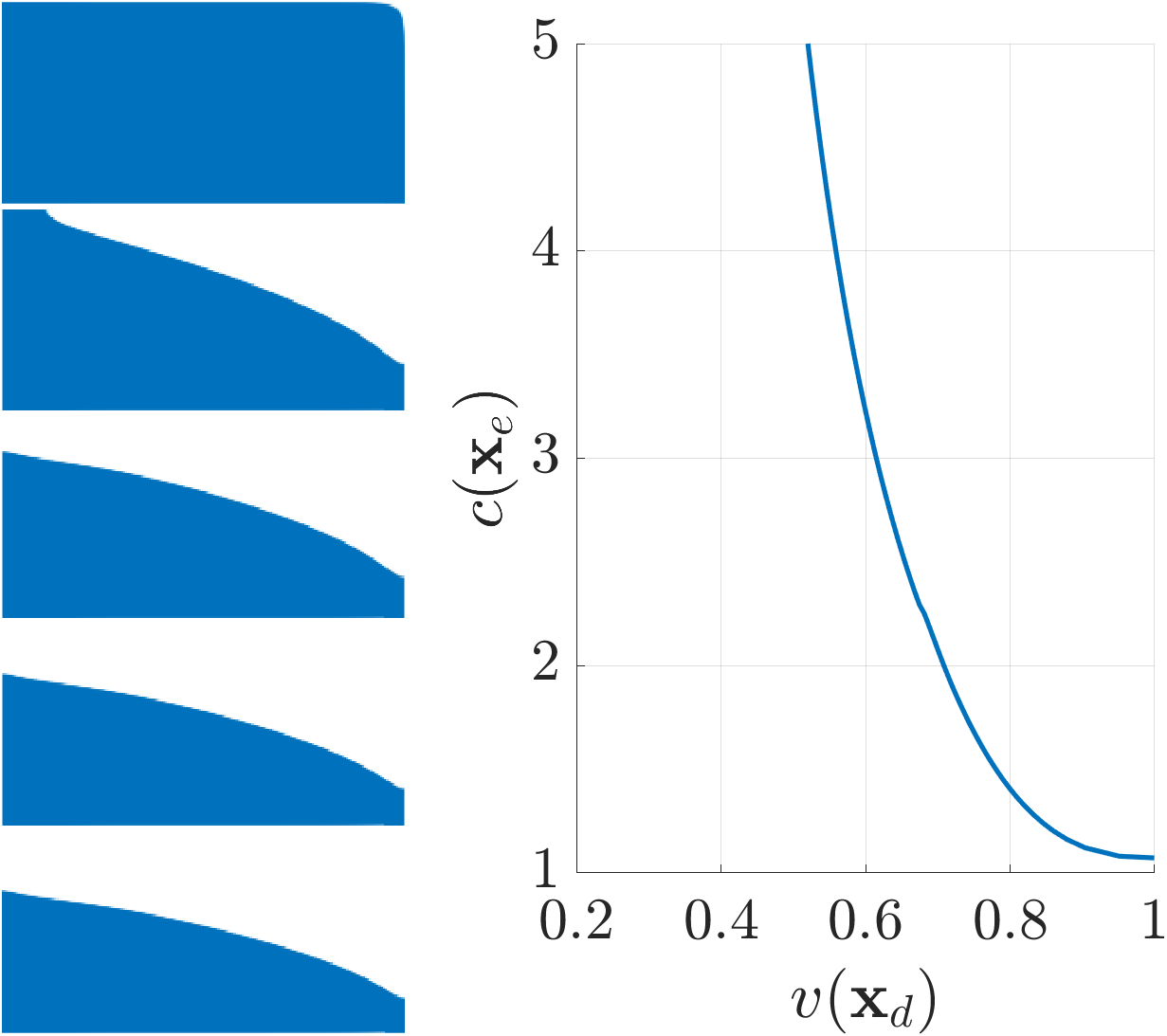}
\caption{Branch 1.}
\end{subfigure}%
\begin{subfigure}[t]{0.33\linewidth}
  \centering
  \includegraphics[width=0.95\textwidth]{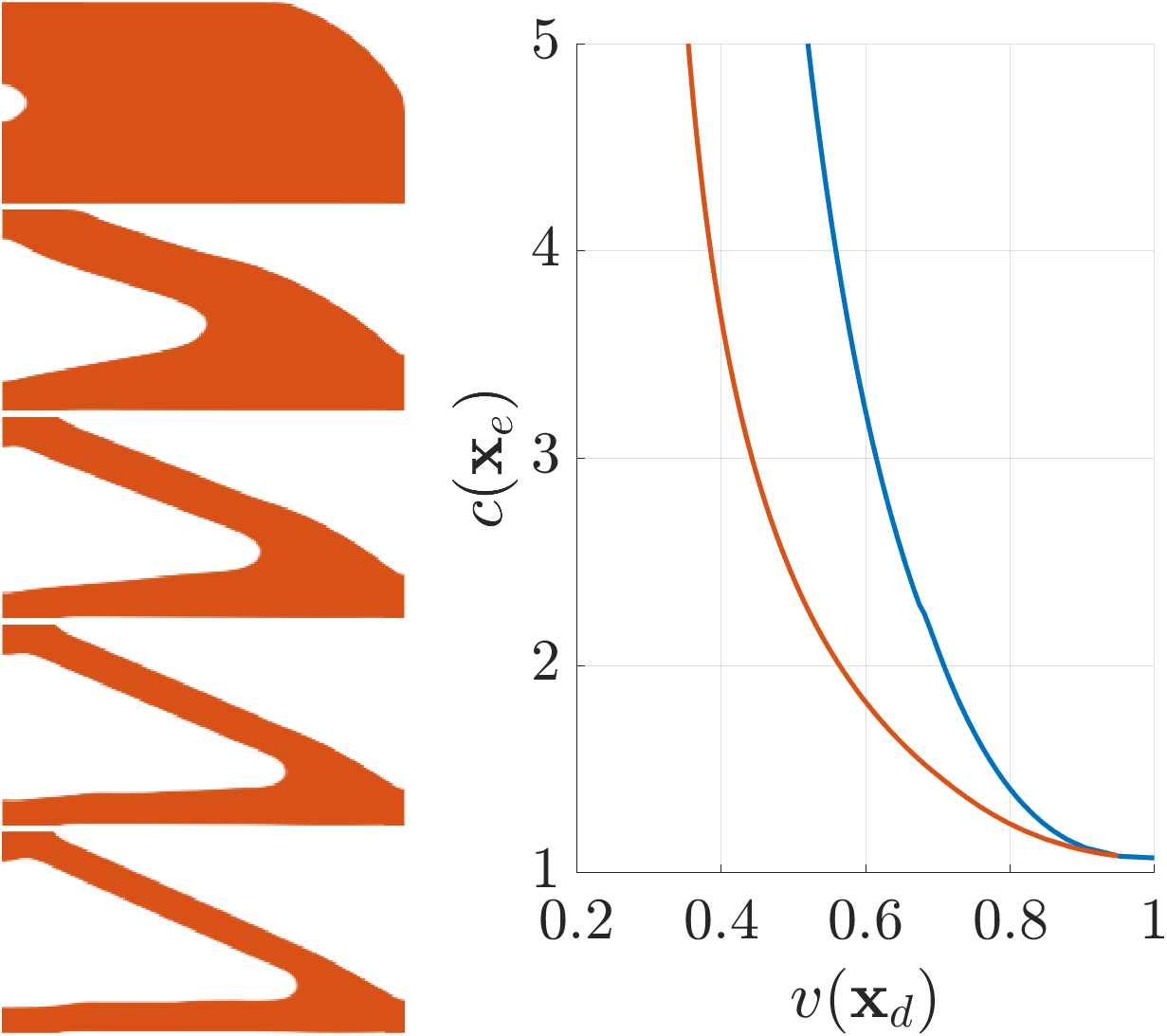}
\caption{Branch 2.}
\end{subfigure}%
\begin{subfigure}[t]{0.33\linewidth}
  \centering
  \includegraphics[width=0.95\textwidth]{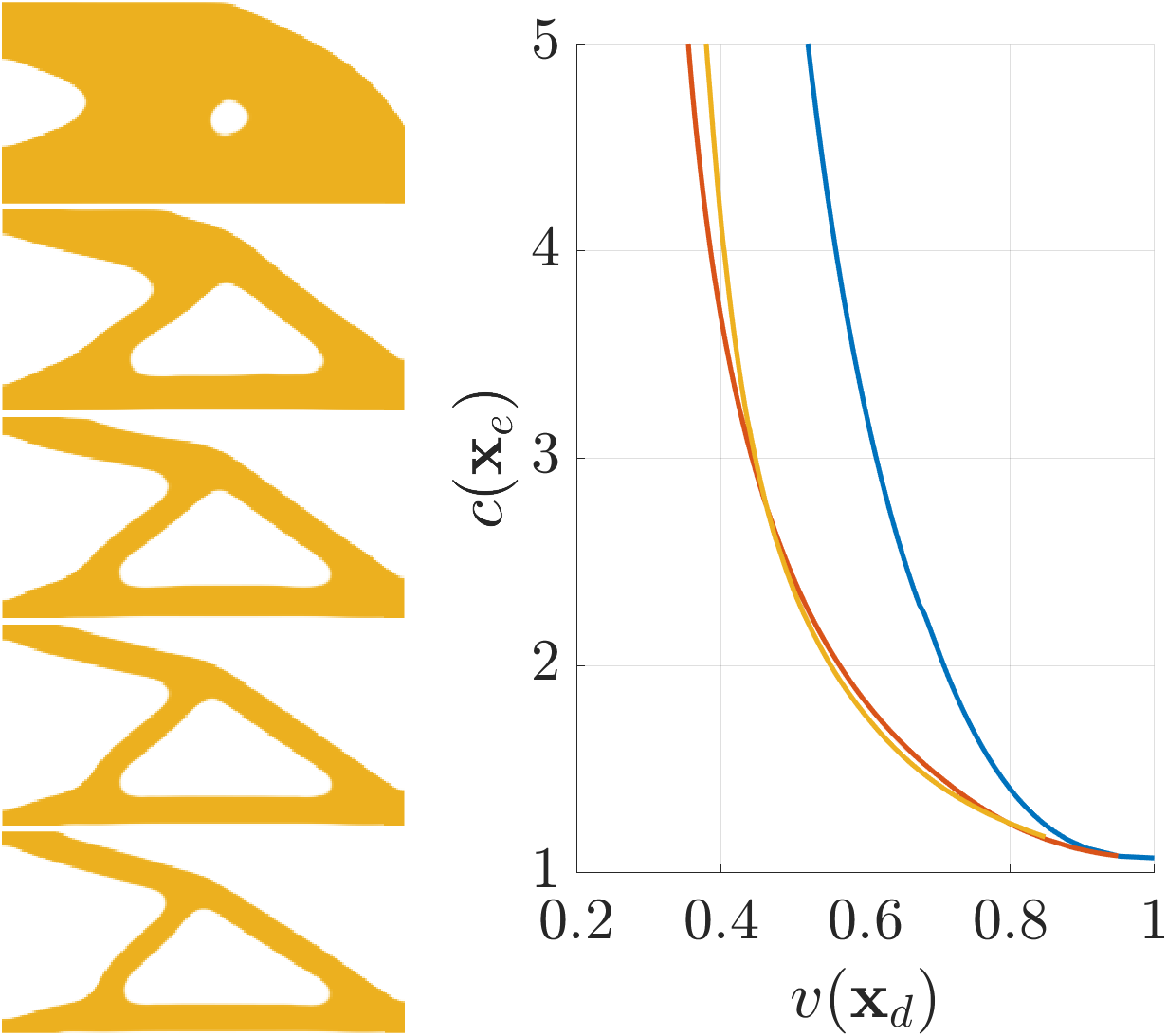}
\caption{Branch 3.}
\end{subfigure}%
\caption{Successive branches for optimization of a $2\times1$ cantilever loaded on the bottom right.\
Each subfigure (a-c) shows a branch, with designs (high to low $v$) on the left and locally nondominated points on the right, in corresponding color.\
Out of the hundreds of obtained locally efficient designs, only four per branch are selected for visualization.\
}
\label{fig:illustration}
\end{figure}

This section studies a $2\,\times\,1$ cantilever beam loaded at the bottom right, to a maximum compliance of $c_{\max}=5$.\
\Cref{fig:illustration} illustrates the branches obtained with the mapping methodology.\
The first branch (blue) takes the fully solid design at $(v,c)=(1,1)$ and brings it to $(0.45, 5)$.\
The design is a single beam connecting the loaded and fixed regions.\
The slight inflection around $(v_d, c_e)=(0.7, 2.2)$ is a boundary effect (see ~\cite{Wallin2020}) that occurs when the beam's edge travels through the top left corner.\\

The second branch (orange) starts from an almost fully solid design nucleated with a hole on the fixed left edge.\ 
When continued to low volume, it transforms into a two-bar truss, ending at $(v,c)\approx (0.4,5)$.\
Except for the original solid design, all nondominated designs now belong to the second branch.\\

The third branch (yellow) has a second hole.\
At low volume, the design resembles a three-bar truss.\
However, for volumes below $v_d=0.45$ and above $0.77$, it is dominated by the two-bar truss.\
Nonetheless, in a wide medium-volume region ($v_d \in [0.45,0.77]$), the three-bar truss dominates, albeat slightly: the gain w.r.t. the two-bar design is a volume reduction of less than $5\%$.\\

Four observations regarding continuity and convexity can now be made.\
First, except for the slight inflection caused by the filtering boundary conditions, all locally nondominated curves are convex.\
Second, the approximation of the nondominated set is \emph{not} convex due to concave kinks when two locally nondominated curves intersect.\
Third, when locally nondominated curves intersect, the nondominated set is nonsmooth, i.e., it has discontinuous derivative.\
Fourth, the nondominated set is a continuous curve: the pieces of locally nondominated curves that make up the nondominated set all cross before terminating.\
Note, however, that this is not guaranteed to hold for setups with less design freedom, as this may cause earlier branch termination.\\

A higher-level trend can also be observed regarding the complexity of the designs, expressed by the number of holes of the efficient designs.\
Namely, the hole count is largest at medium volumes but goes down towards higher and lower volumes.\
For the example in this section, the blue, zero-hole design is nondominated at $v_d=1$.\

\subsection{Study of the high-compliance, low-volume region} \label{sec:highC}

This section investigates the low-volume, high-compliance ($c(\mathbf{x}) \gg 10$) region to verify the observations of the previous section.\
To this end, two cantilevers are studied: one with the force at the bottom right corner, as before, and one with a symmetric setup where the distributed load is applied at the middle of the right edge.\
The maximum compliance $c_{\max}$ is $75$ for the former and $100$ for the latter.\
Only the branches of the zero- and one-hole designs are shown since these are the only ones whose performance is nondominated above $c_e=10$.\\

\begin{figure}[h]
\centering
\begin{subfigure}[t]{0.5\linewidth}
  \centering
  \includegraphics[width=0.95\textwidth]{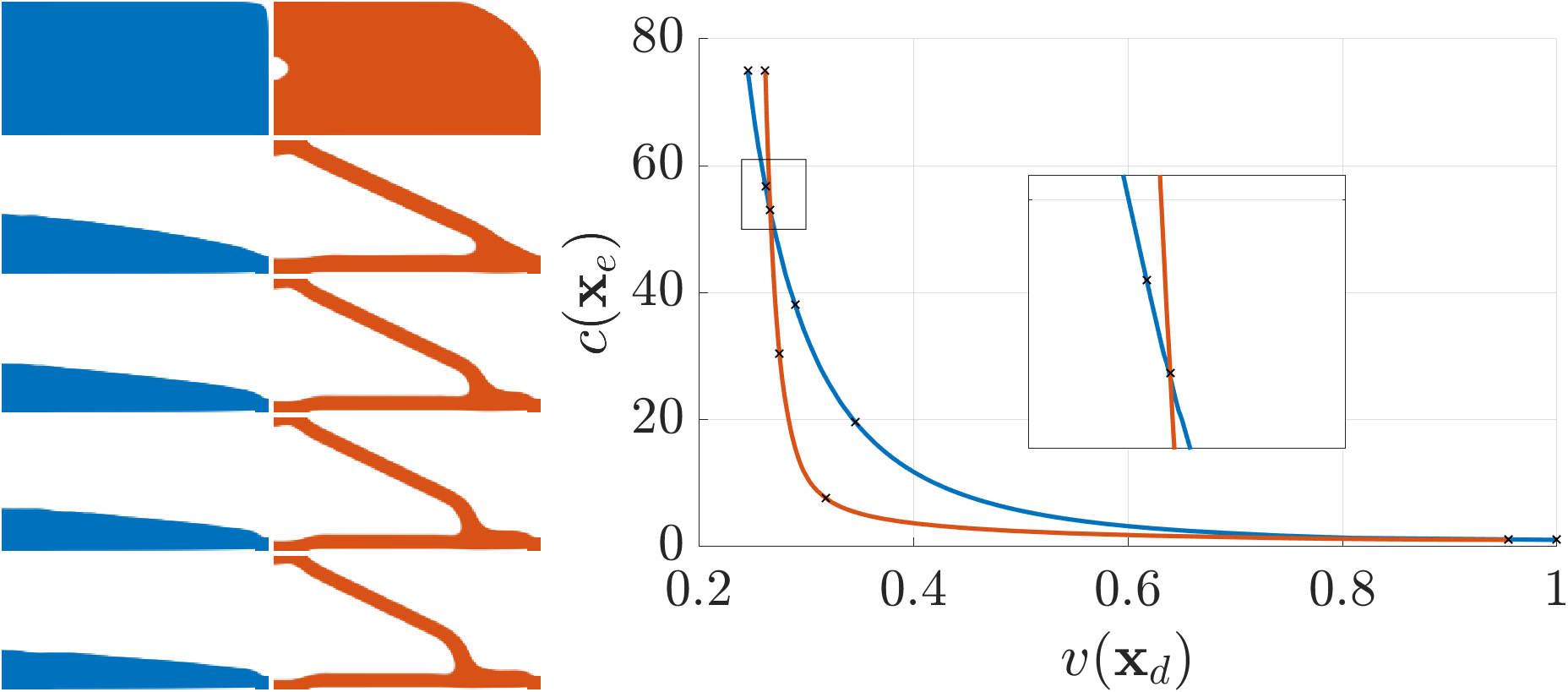}
\caption{$2\,\times\,1$ cantilever, force at bottom right corner.} \label{fig:highC_B}
\end{subfigure}%
\begin{subfigure}[t]{0.5\linewidth}
  \centering
  \includegraphics[width=0.95\textwidth]{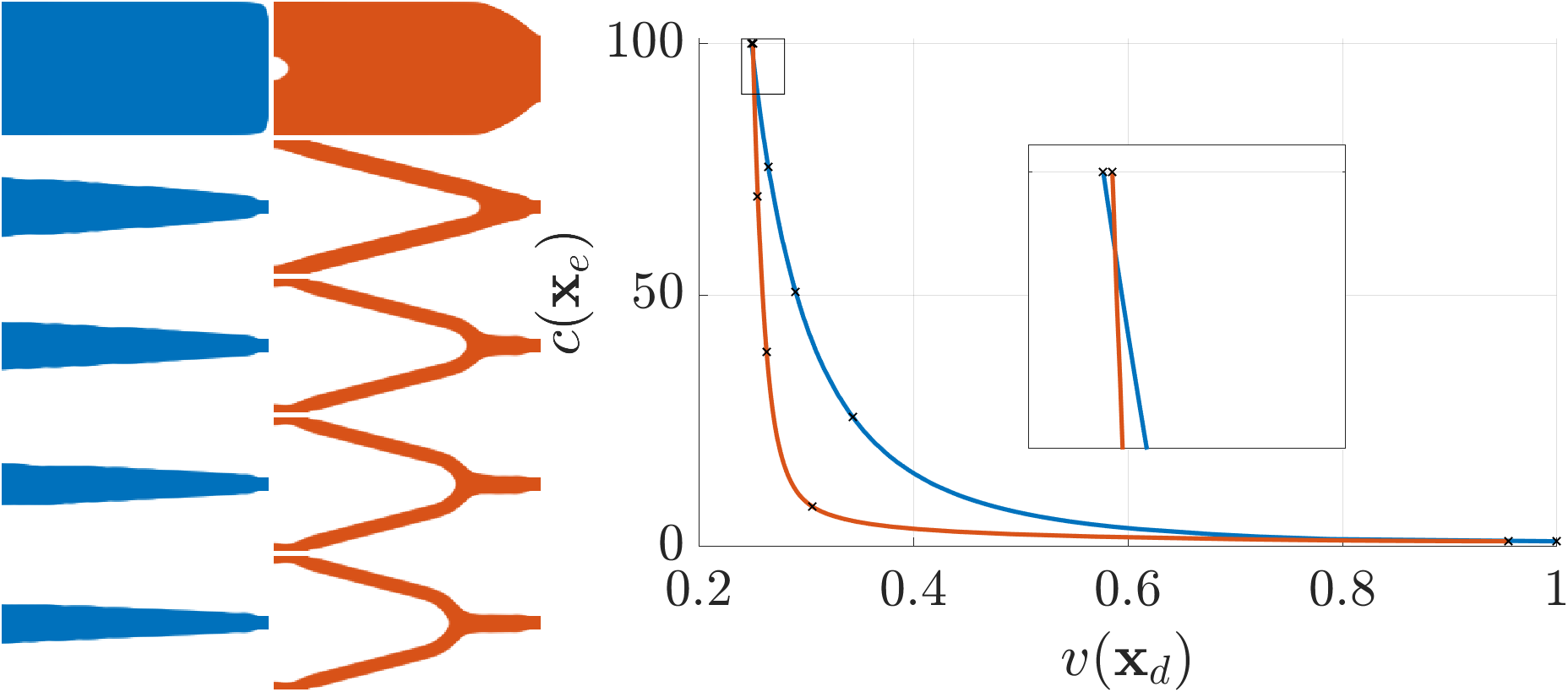}
\caption{$2\,\times\,1$ cantilever, symmetric setup.} \label{fig:highC_M}
\end{subfigure}%
\caption{First two branches (single- and two-bar designs) for a symmetric and asymmetric $2\,\times\,1$ cantilever.}
\label{fig:highC}
\end{figure}

\Cref{fig:highC} shows the obtained results.\
For both the symmetric and asymmetric case, the efficient design switches from the orange two-bar truss to the blue one-bar truss: the blue and orange curves cross at $c_e=53$ and $97$ in \Cref{fig:highC_B,fig:highC_M}, respectively.\
Therefore, the observation that design complexity decreases towards the ends of the nondominated set also holds in the high-compliance region.\\

The convexity observation of the previous section also holds: locally nondominated curves are convex but the nondominated set is not.\
Furthermore, in the high-compliance region the nonconvexity is more outspoken: the blue and orange curves cross at sharper angles than the curves in the medium volume region.\\

The nondominated set also remains continuous: the local frontiers intersect and hence all topological changes along the efficient set occur smoothly in the performance space.\
It must again be noted that such changes \emph{can} occur discontinuously, especially near the outer ends of the nondominated set where the design space is more restricted by the length scale.\
So far, however, any discontinuities are at the order of numerical noise and attributed to the coarse discretization of the locally nondominated curves.\\

\subsection{Symmetric, $2\,\times\,1$ cantilever}\label{sec:symmetry_study_2x1}

 \begin{figure}[h]
     \centering
     \captionsetup{width=\textwidth}
      \includegraphics[width=\textwidth]{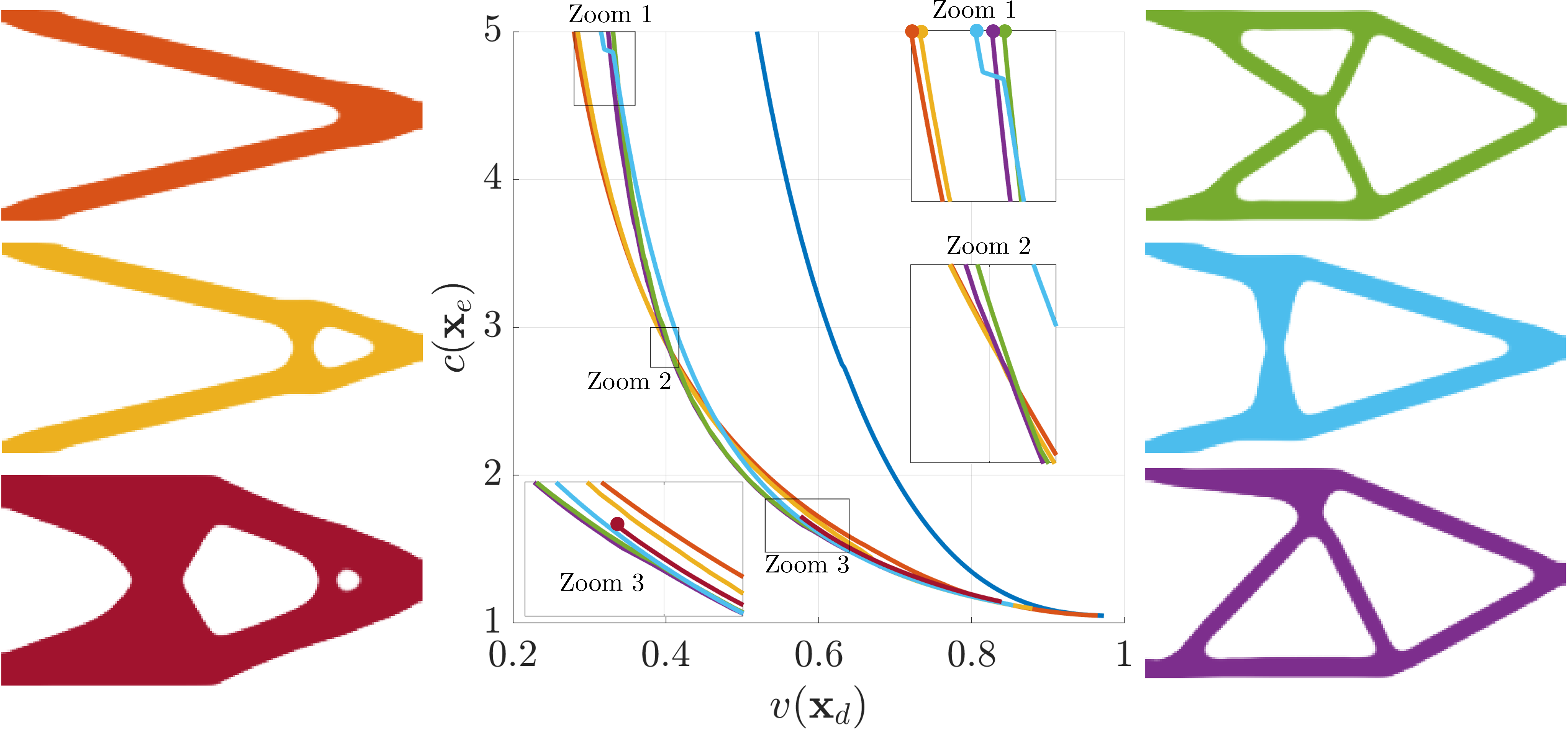}
     \caption{Locally nondominated curves and a selection of designs for a $2\,\times\,1$ cantilever with symmetric boundary conditions.\
     The designs correspond to the dots with the same color in the zoom regions 1 and 3.\
     No designs from the dark blue branch are visualized since these correspond to the 1-bar truss structures of \Cref{fig:highC_M}.}
    \label{fig:canti200}
 \end{figure}

This section considers a $2\,\times\,1$ cantilever under a distributed, symmetric load at the middle of the right edge.\
\Cref{fig:canti200} shows the results of the mapping metholodogy, which led to a total of seven branches.\
The following paragraphs first discuss continuity and convexity and then observe trends regarding the symmetry and complexity of the obtained solutions.\\

The locally nondominated curves are convex and continuous.\
The only exception is (i) the slight inflection of the blue curve, attributed to a boundary effect discussed before, and (ii) the cyan branch's discontinuity near $c(\mathbf{x})=5$.\
This discontinuity is further investigated in \Cref{sec:numerical_examples_2} and linked to a bifurcation.\
As a result, the cyan line essentially connects two locally nondominated curves, each of which is continuous and convex.\\

The approximation of the nondominated set is a continuous curve.\
It consists of segments of the branches, although due to their overlap the exact one is not always clear.\
This is further discussed below.\
Furthermore, note that continuity of the curve is not guaranteed even when the locally nondominated curves are continuous.\
Consider, for example, the case where the mapping methodology only found the branches corresponding to the left-hand designs in \Cref{fig:canti200}.\
The premature termination of the brown branch (see zoom 3) would result in a discontinuity of the nondominated set.\
Finally, note that although the nondominated set is continuous, both its derivative and corresponding designs are not due to changes in the locally nondominated curve.\
For example, in zoom 2 of \Cref{fig:canti200}, the dominant design changes discontinuously from a symmetric five-bar truss to an asymmetric six-bar truss, with a corresponding discontinuity for the slope of the nondominated set.\
As a result, the approximation of the nondominated set is also nonconvex.\\

\begin{figure}[h]
\centering
 \includegraphics[width=0.95\textwidth]{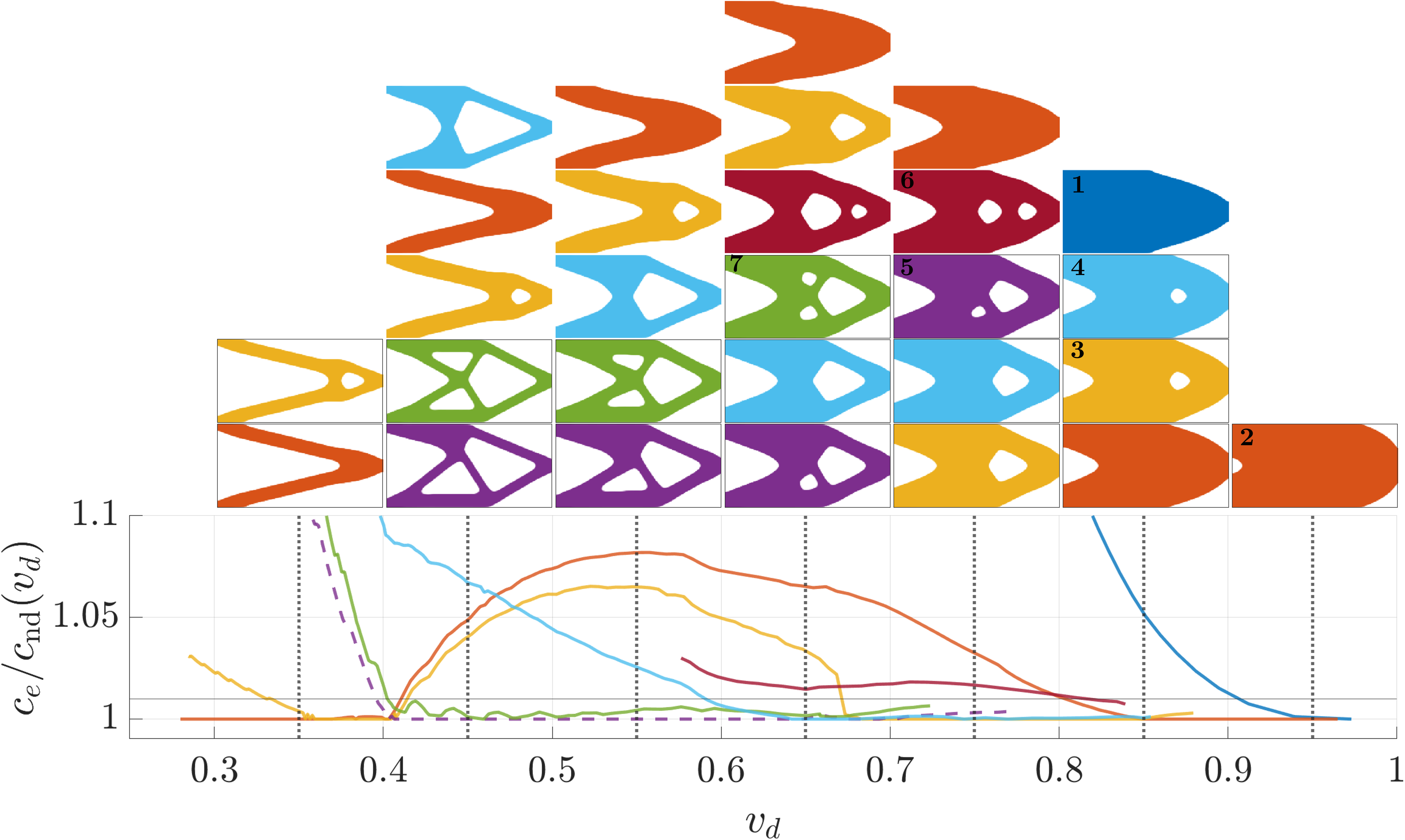}
\caption{Practically nondominated solutions for the $2\,\times\,1$ cantilever. The bottom graph shows the compliance $c_e$ relative to the nondominated design compliance $c_{nd}$, as a function of the volume.\
Only solutions within $10\%$ of $c_{nd}$ are shown.\
Dashed curves correspond to asymmetric designs.\
The top figure shows a selection of the designs at volumes denoted by the vertical dashed lines, ordered by compliance.\
Points with a compliance within $1\%$ of $c_{nd}$ (i.e., below the black horizontal line at $c_e/c_{nd}=1.01$) are considered practically nondominated and the corresponding designs are indicated with a black box.\
For clarity, on top of a color, each branch is also given a number ($1$-$7$) placed in the rightmost plot of the design.
}
\label{fig:symmetry200}
\end{figure}

Many of the curves in \Cref{fig:canti200} are close together.\
This is beneficial for designers if it occurs close to the nondominated set, since it provides multiple options.\
To study this more deeply, \Cref{fig:symmetry200} shows, as a function of volume, the locally nondominated curves that are at most $10\%$ above $c_{\mathrm{nd}}$, the compliance of the nondominated design.\
If the compliance is within $1\%$ ($c_e(v)/c_{nd}(v) < 1.01$), it is considered to be practically nondominated.\
Furthermore, for a sampling of the volumes, the designs are shown and surrounded by a black box if they correspond to a practically nondominated compliance.\\

\Cref{fig:symmetry200} reveals the same trend in design complexity as before: low hole counts at the outer ends of the nondominated set and, conversely, high hole count in the medium volume region.\
For almost all volumes, there are multiple practically nondominated solutions.\\

In recent research~\cite{DeWeer2026}, the symmetry of the solutions for various cantilever lengths was studied.\
There, the pointwise approximation obtained by scalarization prohibited firm conclusions but, for the $2 \times 1$ cantilever, no advantage of asymmetry was observed.\
The mapping methodology of this paper confirms this conclusion but provides a much richer insight.\
Namely, in a wide region $v_d \in [0.4, 0.75]$, both the asymmetric six-bar truss (purple, \# 5) and the symmetric eight-bar truss (green, \#7)  are both practically nondominated.\

\subsection{Symmetric, $4\,\times\,1$ cantilever}\label{sec:symmetry_study_4x1}

\begin{figure}[h]
     \centering
     \captionsetup{width=\textwidth}
      \includegraphics[width=\textwidth]{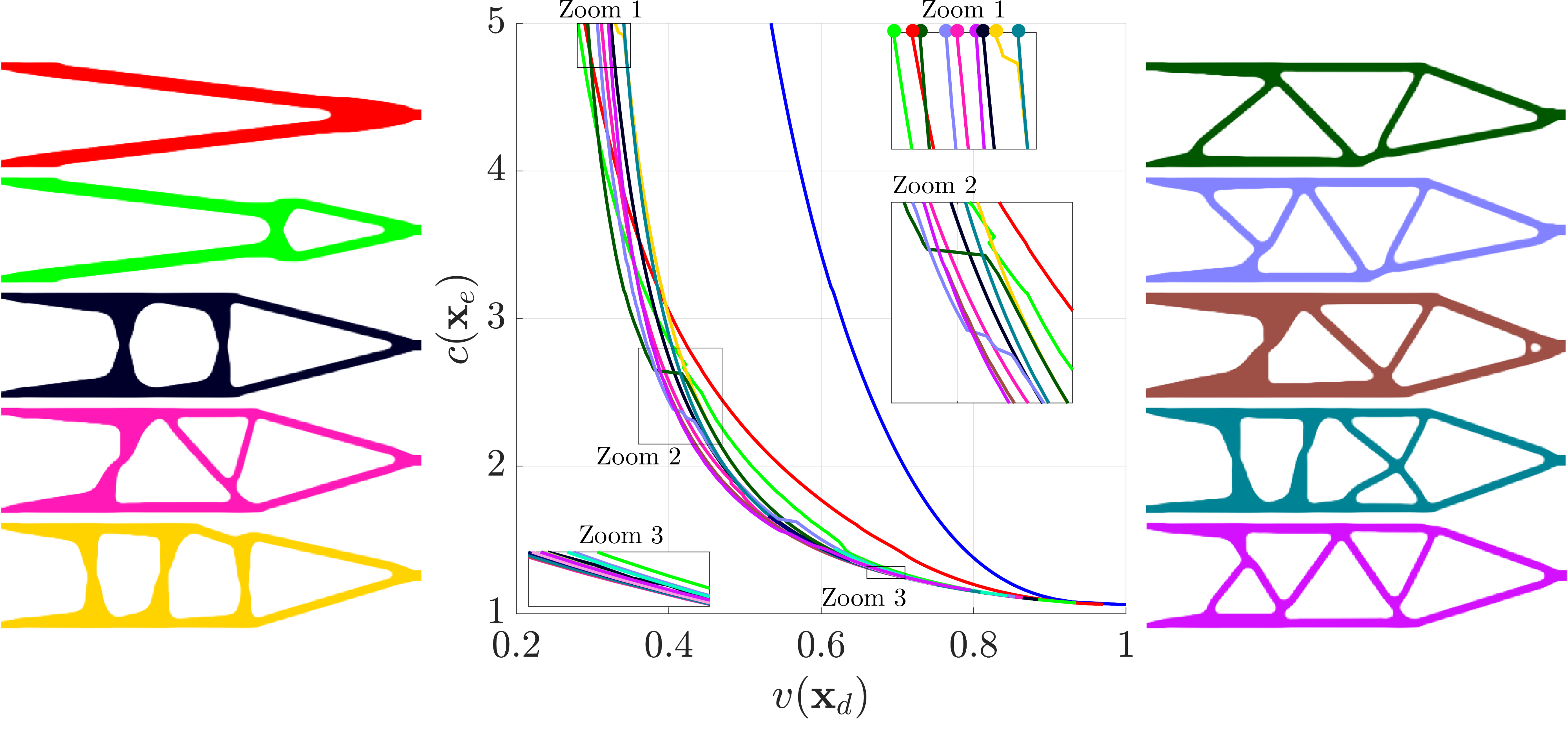}
     \caption{Locally nondominated curves for a $4\,\times\,1$ cantilever with symmetric boundary conditions.\
 A total of fifteen branches are visualized, each with a different color.\ 
Eight designs are visualization to the left and right of the frontier plot.\
They correspond to the dots with the same color in zoom 1 and all attain the maximum compliance bound.}
    \label{fig:canti400}
 \end{figure}

The previous section confirms that, for short cantilever lengths, the benefit of asymmetric designs is negligible.\
However, recent research also indicated that asymmetry can outperform symmetric designs for longer cantilevers.\
To confirm this, the current section applies the branch mapping methodology to a $4\,\times\,1$ cantilever, with results in \Cref{fig:canti400}.\\

 The $4\,\times\,1$ cantilever shows broadly the same behaviour as the $2 \,\times\,1$ cantilever: (i) except for bifurcations, the locally nondominated curves are convex and continuous , and (ii) the nondominated set is a continuous but nonsmooth curve.\
 However, there are two significant differences.\
First, branch terminations occurs more often: of the fifteen branches, only nine reach the maximum compliance bound.\
Termination occurs either (a) due to collisions with an existing branch or (b) due to a topological change.\
Second, many more bifurcations are noticeable: zoom region 2 shows several successive bifurcations.\
Bifurcations are studied more deeply in \Cref{sec:numerical_examples_2}.\
We here simply note that the slight discontinuities of the nondominated set are caused by these numerical artifacts but likely do not reflect its true behaviour.\\

\begin{figure}[h]
\centering
 \includegraphics[width=0.95\textwidth]{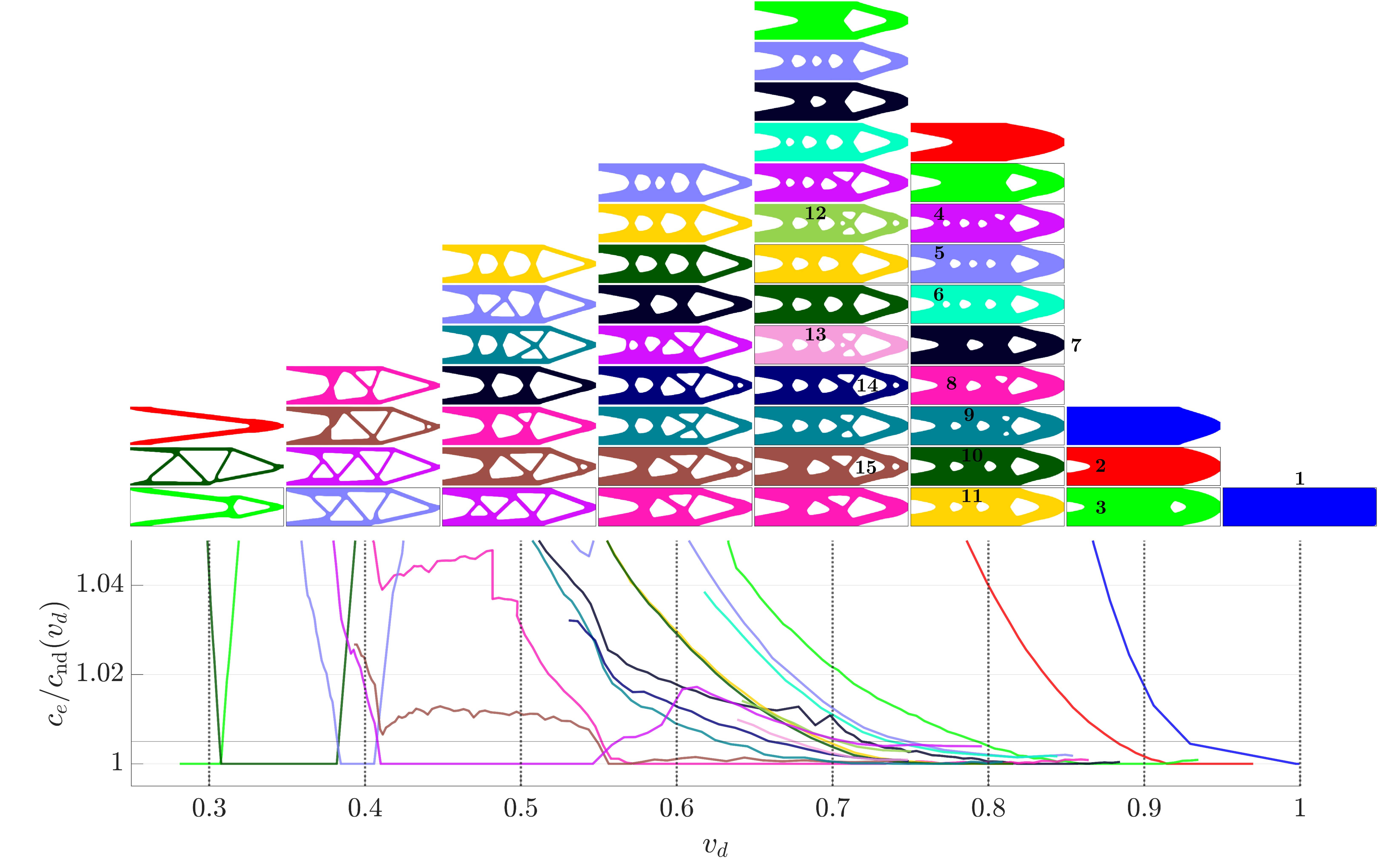}
\caption{Practically nondominated solutions for the $4\,\times\,1$ cantilever. The plot mirrors \Cref{fig:symmetry200} except for a different tolerance for practical nondominance ($0.5\%$ instead of $1\%$).\ 
Only designs within $5\%$ of $c_{nd}$ are shown.\
}
\label{fig:symmetry400}
\end{figure}

\Cref{fig:symmetry400} studies the practically nondominated points for the $4\,\times\,1$ cantilever, which now includes points within $0.5\%$ of $c_{\mathrm{nd}}$.\
Compared to the $2 \times 1$ cantilever, the pattern is broadly similar: below $v_d=0.55$, only one design is practically nondominated.\
This briefly increases to two designs in $v_d \in [0.55,0.65]$, after which all branches start colliding and the number of options increases dramatically to seven designs at $v_d=0.7$.\
The, densification occurs and branches start terminating until only the completely solid design remains.\\

A notable difference is noticeable regarding the symmetry of the designs corresponding to practically nondominated points.\
Whereas the $2\,\times\,1$ cantilever always offered at least one competitive symmetric design, this is no longer the case for the $4\,\times\,1$ cantilever.\
Between symmetric design \#3 at $v_d\approx 0.3$ and symmetric design \#9 at $v_d\approx0.65$, there are no practically nondominated points corresponding to symmetric designs.\
Instead, the $v_d \in [0.3,0.65]$ region is dominated, in succession, by (i) an asymmetric 9-bar truss (design \#10), (b) an asymmetric 11-bar truss (design \#5), (c) an asymmetric 13-bar truss (design \#4) and then (d) the asymmetric designs \#8 and \#15.\
In short, there is a wide region where asymmetric designs outperform symmetric ones.\

\section{Numerical examples: investigating instabilities caused by flatness and bifurcations} \label{sec:numerical_examples_2}

The previous section showed branches for several setups, which revealed two types of numerical instabilities.\
First, some trees show branches with a very similar layout, which signifies a lack of uniqueness.\
Second, some branches exhibit discontinuities, which stems from a sudden loss of optimality.\
This section studies these numerical instabilities in more detail.\

\subsection{Flatness} \label{sec:flatness}

So far, a clearly observed phenomenon is the appearance of designs with the same topology but slightly different layout.\
Their bars and holes are located in the same general location, but moved one or more element widths -- high enough to not trigger the equality check detailed in \ref{app:continuation}.\\
 
An illustrative example is the symmetric $2\,\times\,1$ cantilever study, summarized in \Cref{fig:symmetry200}.\
Therein, the branch mapping methodology obtains two variants of a five-bar truss: the yellow (\#3) and cyan  (\#4) design at $v_d=0.75$.\
Although the designs are visually similar, they are numerically distinct: the vertical bar's horizontal location differs by several element widths between the two designs.\
Nonetheless, their performance is very close: around $v_d=0.75$ the yellow and cyan lines are indistinguishable.\
This signifies a local flatness of the optimization landscape and raises questions about the uniqueness of an optimum.\
Furthermore, the two branches diverge significantly towards lower volumes: the yellow branch forms a two-bar truss punctured near the load, whereas the cyan branch forms a five-bar truss.\
This signifies a strong dependence on the initial conditions.\
Because this has repercussions for branch mapping, they are analyzed further in this section.\\

\begin{figure}[h]
\centering
  \includegraphics[width=0.95\textwidth]{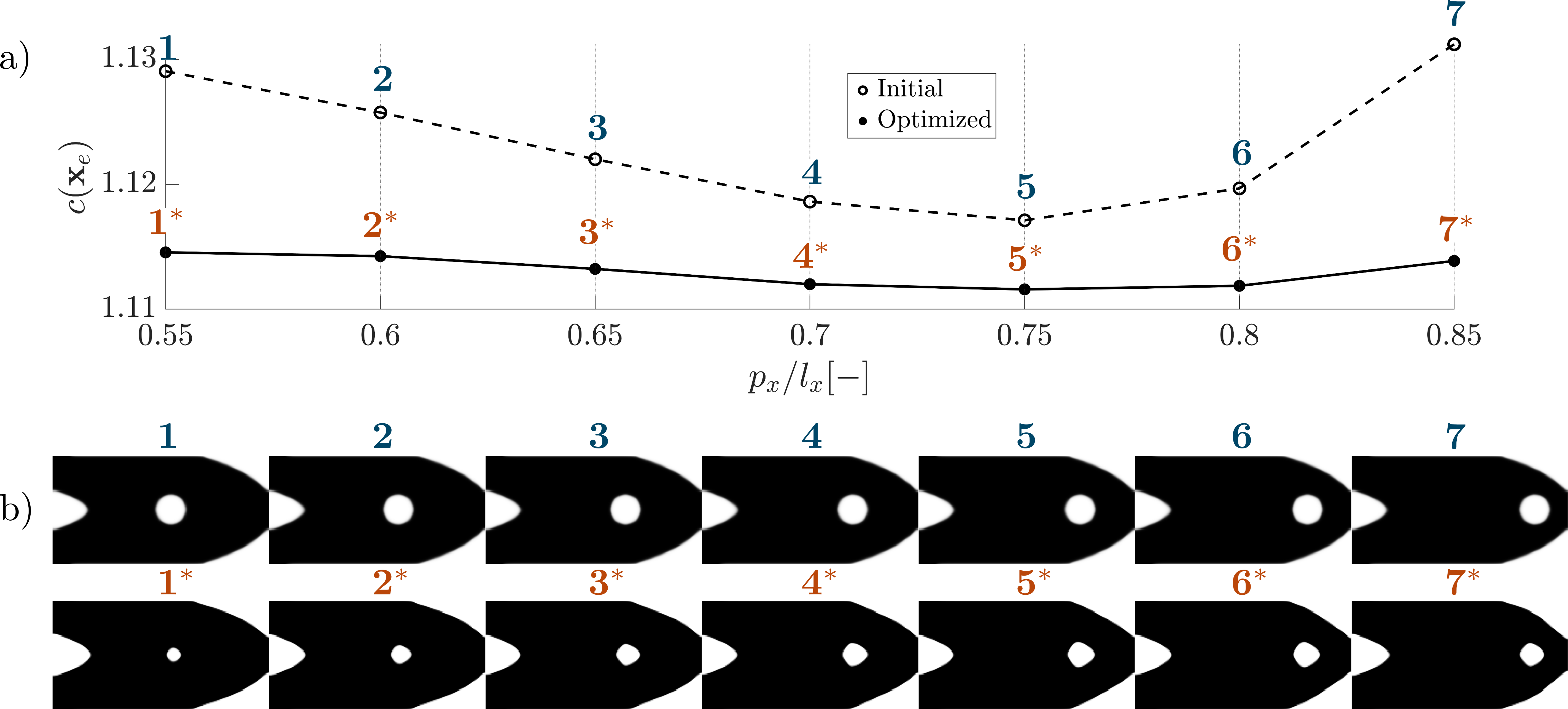}
\caption{Result of seven compliance minimization runs, labelled $1-7$. Each run uses the same maximum volume bound ($v(\mathbf{x}_d) \leq v_{\max}=0.865$) and starts from an initial guess with a hole along the symmetry axis, whose horizontal position $p_x$ varies. (a) Compliance (initial and optimized) as function of $p_x$. (b) Initial (top) and optimized designs (bottom). }
\label{fig:flatness}
\end{figure}

To illustrate local flatness, a two-bar design is punctured with a circular hole at $(p_x, 0)$, i.e., along the symmetry axis.\
For seven values of $p_x/l_x$, a compliance minimization with the same maximum volume $\varepsilon_v=0.865$ is performed.\
\Cref{fig:flatness} shows the result of these runs, labelled $1$ to $7$.\
All optimized designs retain the same topology by adapting the shape of the hole while keeping its center at roughly the same location.\
The optimized topologies have reduced compliance, with the best performance achieved by the fifth run (denoted $5^*$).\
Nonetheless, compared to the initial guesses, the optimized compliance is very similar for all runs.\
Whereas the initial compliances lie within a $0.01$ range, the optimized compliances are less than $0.003$ apart.\\

All optimized designs except $5^*$ could be further improved by moving closer to design $5^*$.\
The absence of this behaviour can be attributed, at least in part, to the local flatness of the objective landscape.\
From a methodological perspective, however, it may also be attributable to the optimizer itself.\
Whereas the MMA optimizer allows small increases in the objective during optimization, the GCMMA optimizer employed in this work is more conservative.\
Iterations are only accepted when they decrease the objective.\
This benefits the branch tracking, since it ensures local minima are followed closely throughout the multiobjective optimization landscape.\
However, the smooth compliance curves in \Cref{fig:flatness} may still contain a significant amount of numerical noise.\
In that case, small local irregularities would prevent GCMMA from converging to design $5^*$.\\

Several future research directions could be explored to deal with local flatness.\
First, a methodological improvement could be to tune the numerical relaxation, e.g. the $\beta$ projection, to balance accurate branch tracking against the ability to traverse small local irregularities.\
Second, to avoid repeating branches, the branch-collision tolerances could be adapted in function of the local curvature: tight tolerances in regions of high curvature and looser tolerances in locally flat regions.\
This requires computing the local curvature, i.e., the Hessian, which is computationally tractable via a Newton-based approach with simultaneous analysis and design~\cite{RojasLabanda2015}.\
Third, instead of employing branch-collision tolerances, one could instead turn to deflation to avoid local minima, as done by Papadopoulos et al.~\cite{Papadopoulos2021}.\
Finally, since the local irregularities might be attributed to the material interpolation, a fourth research direction could focus on employing a conformal mesh to avoid the interpolation altogether~\cite{Allaire2014}.\\

\subsection{Bifurcations}\label{sec:bifurcations}

\Cref{def:local_frontier} requires that, in order for two locally nondominated points to be part of a locally nondominated curve, they should be ``connected''.\
Specifically, there should be a transformation between them that preserves topology and local optimality.\
This provides two failure points for a branch's continuation.\
First, if the design changes topology during continuation, then a transformation would violate the length scale.\
The topology is checked numerically and causes most of the branch terminations.\
Second, a design can retain its topology but temporarily lose local optimality.\
The result is a local optimum with the same topology but strongly different performance.\
Numerically, this is noticed as an inflection/discontinuity of the locally nondominated curve.\
This is now further investigated, as it signifies a bifurcation and the possibility of branch-splitting.\\

\begin{figure}[h]
\centering
  \includegraphics[width=0.95\textwidth]{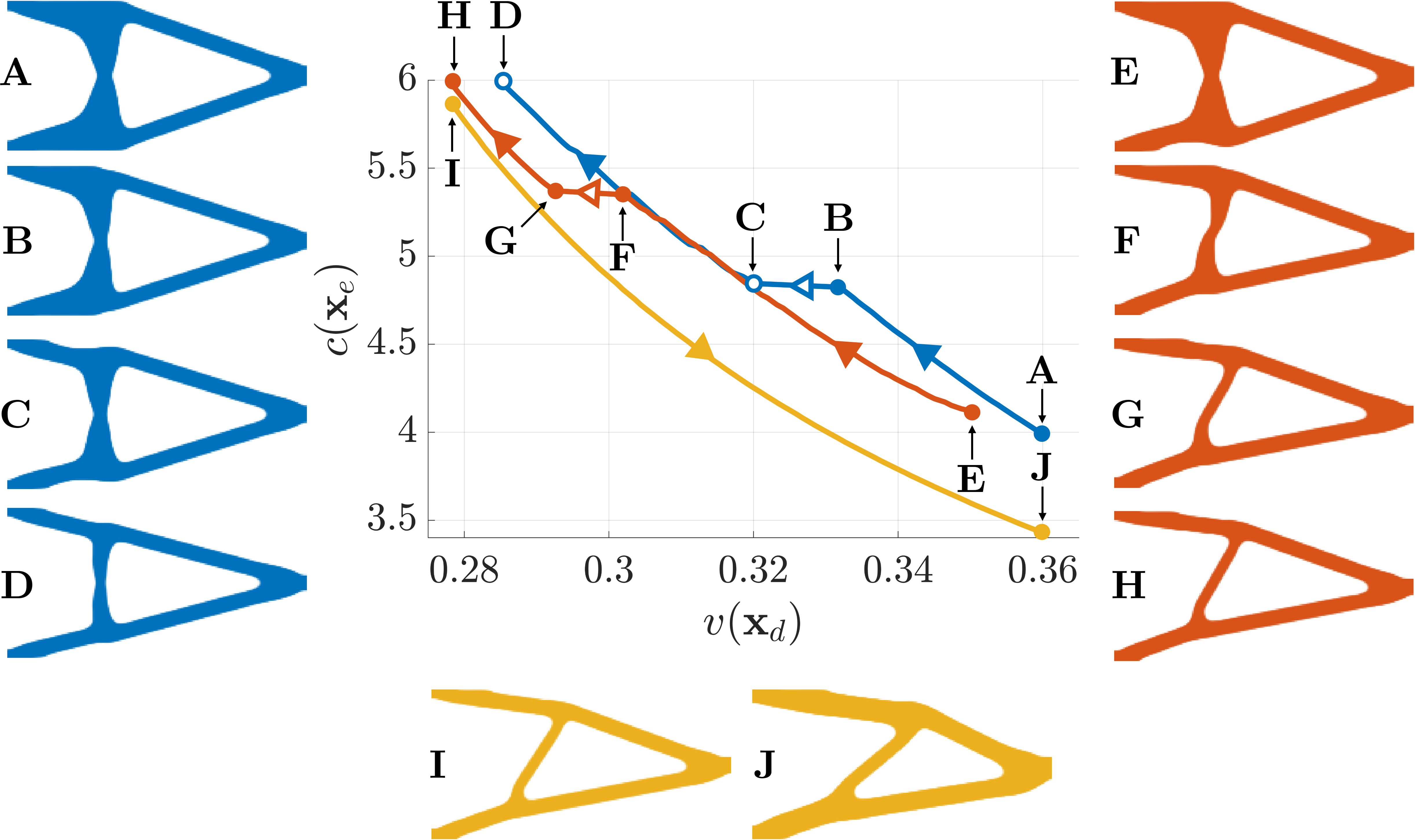}
\caption{}
\label{fig:fivebar_study}
\end{figure}

To frame the phenomenon studied in this section, the different types of branch terminations observed so far are now categorized:
\begin{enumerate}
\item A continuation step leads to a design with a \emph{different} topology:
\begin{enumerate}
\item An increase of the maximum volume $\varepsilon_v$ fills a hole.\
The genus drops and continuation terminates.\
\item An increase of the maximum compliance $\varepsilon_c$ removes a bar that (usually) was at the minimum length scale.\
The genus drops and continuation terminates.\
\end{enumerate}
\item A continuation step leads to a design with the \emph{same} topology:
\begin{enumerate}
\item A change in $\varepsilon_c$ or $\varepsilon_v$ moves a design's boundary along a changing filtering boundary condition.\
The locally nondominated curve inflects slightly.\\
\textit{Example: blue branch in} \Cref{fig:illustration}.\
\item An increase of $\varepsilon_c$ compresses the design's joints together, shortening the connected bars.\
The volume reduces strongly and leads to a discontinuity in the locally nondominated curve.\\
\textit{Example: the cyan five-bar truss in} \Cref{fig:symmetry200}.\
\end{enumerate}
\end{enumerate}
This section studies the last category, option 2b, by investigating the cyan five-bar truss in \Cref{fig:symmetry200}, which undergoes an inflection around $c_e=5$.\
\Cref{fig:fivebar_study} shows the obtained designs and curves in the $(v_d,c_e)$ space.\\

Design A is used as starting point for the first continuation: it is a five-bar truss corresponding to the cyan branch of \Cref{fig:symmetry200}, at $c_e=4$.\
It is continuated to $c_e=6$ via a hundred successive volume minimization steps.\
At $c_e=4.82$, a bifurcation occurs: the joints connecting the inner bar of design B disconnect from the outer boundary.\
This straightens the outer bars and shortens the inner bar, leading to a strongly reduced volume.\
Further continuation accelerates this process until design D is reached.\
However, since local optimality is temporarily lost between designs B and C, the result is effectively two separate subbranches.\\

A further investigation of the top subbranch, section C-D, reveals that the optima are highly unstable.\
That is, starting a compliance minimization run with C or D as initial guess and using as $\varepsilon_v$ the corresponding volume does not lead to C or D, respectively.\
Instead, the inner bar is immediately removed and a two-bar design is obtained.\
The designs along the C-D section might thus be closer to saddle points than local minima.\\

The bifurcation between C and D seems caused by the symmetric boundary conditions: the design sticks to the top and bottom boundaries until it is no longer possible, after which it disconnects from both boundaries simultaneously.\
To investigate this further, a second branch is created by starting from a particular initial guess: design A with a semicircular hole at the bottom edge.\
Volume minimization then yields designs E to H.\
Due to the asymmetric initial guess, design E no longer sticks to the bottom edge.\
As a result, the eventual bifurcation from F to G, where it disconnects from the top edge, occurs later (at $c_e=5.35$).\
Furthermore, the top section G-H is now stable.\

A third continuation run is now performed by starting from design H and minimizing compliance.\
No bifurcations are observed, so the section I-F is a single branch.\
Furthermore, although it seems that the symmetric A-B branch and the asymmetric I-F branch will eventually cross in the $(v_d,c_e)$ space, their designs are fundamentally different and hence will not cross.\
Finally, the difference between designs H and I is significant: their volume is the same (up to optimization tolerance) but the compliance is not ($c_e=6 \rightarrow 5.89$).\
This is due to the inner bar, which slightly rotates clockwise.\
Thus, although volume minimization with a compliance constraint (\Cref{eq:minV_variant}) and compliance minimization with a volume constraint (\Cref{eq:minC_variant}) are made mathematically equivalent up to a choice of $\varepsilon_c$ and $\varepsilon_v$, their numerical behaviour deviates.\\

Taken together, the numerical experiments presented in this section show several ways in which volume minimization can lead to successive bifurcations that break symmetry.\
Namely, the transitions B-C and F-G see an asymmetric change of the design: a disconnection from the design domain boundary and a rotation of an internal bar.\
It is important to note that these symmetry-breaking design transitions could also be applied in the other direction: designs E to H could be flipped along the symmetry axis without changing neither their optimality nor their performance.\
This discontinuous symmetry-breaking behaviour is reminiscent of subcritical pitchfork bifurcations.\\

It is conceivable that the inverse can also occur: a transition from asymmetry to symmetry as the compliance decreases.\
Consider for example the extension of points A and J to lower compliance.\
Assuming neither curve terminates prematurely, the corresponding blue and yellow curves will either merge or cross.\
The former indicates a smooth transition from symmetry to asymmetry, reminiscent of a supercritical pitchfork bifurcation.\
The latter will again lead to a subcritical pitchfork bifurcation, where the asymmetric design loses optimality and converges to the dominant symmetric design.\
The result is then a hysteresis between symmetry and asymmetry.\
Instead of further studying this phenomenon, we here conclude this study by noting that bifurcations can indeed occur and that their exact behaviour depends on parameters such as length scale, topology and cantilever length.\\

Although bifurcations are physically significant phenomena, the mere possibility of their occurence is also important from a methodological perspective.\
The simple methodology of this paper, for example, assumes that topological connectivity is enough to ensure that local optimality is also conserved.\
Consequently, locally nondominated curves are drawn with inflections and discontinuities even though these are now understood to indicate a loss of optimality.\
A more correct visualization of these curves would be to draw them separately, perhaps only connected with a dashed line to indicate their topological similarity.\\

We now outline two similar but opposed directions for future research of bifurcations in multiobjective topology optimization.\
The first direction is to consider bifurcations detrimental to numerical stability, since they strongly increase the iteration count of a continuation step.\
Bifurcations should thus be prevented by detecting them before they occur.\
A good indicator of an impending bifurcation is the local curvature.\
Reconsider for example the classic pitchfork bifurcations in \Cref{fig:pitchforks}: as $\lambda \underset{<}{\rightarrow} 0$, the second order terms of $f(x)$ start cancelling until only higher-order terms remain at $\lambda=0$.\
Bifurcations can thus be avoided by tracking local curvature, leading to similar suggestions as in \Cref{sec:flatness}.\
The second research direction starts from the opposite conclusion: instead of treating bifurcations as undesirable, they are considered beneficial.\
After all, despite the numerical instability and accordingly higher iteration counts, bifurcations lead to a significant performance improvement.\
Therefore, intentionally inducing bifurcations could be a good way to avoid poorly performing local optima, as illustrated in recent works on other applications~\cite{DeWeer2025,Elbek2026}.\

\section{Conclusions} \label{sec:conclusions}

This paper studies multiobjective topology optimization.\
A previous study showed that the nondominated set, also referred to as the Pareto frontier, consists of segments of locally nondominated curves~\cite{DeWeer2026}.\
In contrast to the continuous and convex Pareto frontiers often reported in the literature, this piecewise structure can give rise to nonconvexities and discontinuities in the nondominated set.\
Although scalarization methods give rise to pointwise approximations that often obscure these properties, the numerical examples in \cite{DeWeer2026} show nondominated sets that exhibit signs of nonconvexities, discontinuities and an underlying branch-like structure.\ 
The goal of this paper is to further uncover this structure by mapping its underlying branches.\
To limit the scope, it is focused on minimization of volume and compliance.\\

The study provides three main contributions: a methodology to map branches and three physical observations related to properties of the nondominated set:
\begin{itemize}
\item The methodological contribution consists of a simple algorithm to map branches.\
It contains two parts: (i) a hole nucleation heuristic, based on topological derivatives, to quickly generate promising new designs, and (ii) a continuation of the $\varepsilon$-constraint scalarization parameters.
\item The first physical observation, which confirms and expands the findings in \cite{DeWeer2026}, relates to convexity and continuity.\
Namely, the locally nondominated curves are all convex, continuous and smooth.\
The nondominated set, which is built from segments of locally nondominated curves, remains continuous but can be nonsmooth and nonconvex when these curves intersect.\
At these intersections, the change in the design space is also discontinuous.\
Discontinuities are not observed but are conjectured to occur for cases with lower design freedom, as this can cause premature branch termination.\
 \item The second physical observation relates to the symmetry and complexity of the solutions.\
 Due to overlapping branches, especially at medium-high volumes, many designs have practically the same performance.\
 This provides designers with several suitable options.\
 However, near the outer edges of the frontier, the length scale causes branch termination, which culls the number of options.\
 Furthermore, for some regions asymmetric designs are observed to outperform symmetric designs even when the boundary conditions are symmetric.\
 This confirms the observation in \cite{DeWeer2026} and the references cited therein.\
 \item The third and final observation concerns two properties of the multiobjective optimization landscape.\
 First, the numerical examples reveal a recurrence of similar designs.\
 This lack of uniqueness is linked to an observed flatness in the objective landscape: large topology-preserving design changes can result in only minimal objective variations.\
Second, the numerical examples reveal several discontinuities in the locally nondominated curves.\
These do not reflect any true discontinuity but instead arise from bifurcations.\
Specifically, these bifurcations are associated to symmetry-breaking subcritical pitchfork bifurcations.\
As far as the authors are aware, this work is the first to demonstrate branch-splitting behaviour in multiobjective topology optimization.\
\end{itemize}

Several future research directions could be explored to extend the conclusions of this paper.\
First, it would be interesting to extend the approach towards different objectives (maximum stress, buckling, ...) because these settings lead to locally nondominated curves with different characteristics~\cite{DeWeer2026}.\
Second, an extension to three or more objectives is also valuable but methodologically challenging, since it requires continuation in additional dimensions.\
Third and final, both extensions would benefit from a study of the observed numerical instabilities.\
For example, a combination of deflation~\cite{Papadopoulos2021} and higher-order continuation~\cite{Cesarano2025,Gangl2026} could be used to monitor uniqueness, avoid previous solutions and more tightly control the optimizer path.\
Additionally, a conformal mesh~\cite{Allaire2014} and/or advanced continuation strategies~\cite{DeWeer2025,Elbek2026} could be employed to circumvent local irregularities in the multiobjective optimization landscape.\

\section*{Statements and Declarations}
\subsection*{Funding}
Resources and services were provided by the VSC (Flemish Supercomputer Center), funded by the Research Foundation - Flanders (FWO) and the Flemish Government.\

\subsection*{Conflict of Interest}

On behalf of all authors, the corresponding author states that there is no conflict of interest.

\subsection*{Author Contributions}
All authors contributed to the study conceptualization and design, commented on previous versions and read and approved the final version.\
Data collection, analysis and writing was performed by the corresponding author.\

\subsection*{Ethics approval and Consent to participate}

Not applicable.

\subsection*{Data Availability and Replication of Results}

Methodological details are described to facilitate the replication of the results.\
Data and/or Matlab codes will be made available upon reasonable request.

\appendix

\section{Topology optimization} \label{app:TO}

This appendix lists numerical details regarding the topology optimization routine.\
Parameters that vary between the numerical examples are listed in \ref{app:values}.\\

The optimization framework is written in Matlab and run on the Flemish Supercomputer.\
During continuation, the same optimization routine is called with different parameters (initial guess, $\varepsilon$ values, ...).\
The material is interpolated with SIMP~\cite{Bendsoe1989}:
\begin{equation}
x \rightarrow \delta + (1-\delta)x^p,
\end{equation} 
where the penalization constant is $p=3$ and $\delta=1e-6$.\
All experiments use isotropic material with Young's modulus of $1\,\mathrm{N}/\mathrm{m}^2$ and Poisson ratio of $0.3$.\\

As discussed, the design variables are filtered with a PDE filter with radius $R$, rescaled with a $\sqrt{3}$ factor to ensure it corresponds to a classic convolutional hat filter~\cite{Lazarov2011}.\
Then, projection with steepness $\beta$ ensures a crisp design: if $\beta=2R/h$, with $h$ the element width, then a single gray element layer remains between the solid and void regions.\
For numerical stability, however, $\beta=10$ for all examples and $R$ varies s.t. $\beta=R/3h$ or $\beta=R/5h$ (see \ref{app:values}).\
This aids numerical stability and retains the length scale, but at the cost of several (two to three) layers of gray elements between void and solid regions.\\

Locally nondominated points are found by solving a scalarized version of \Cref{eq:min_cv} to local optimality.\
The compliance minimization variant is implemented as stated in \Cref{eq:minC_variant}.\
The volume minimization variant uses a normalization of the compliance constraint:

\begin{equation} \label{eq:minV_variant_actual}
  \min_{\mathbf{x}} v(\mathbf{x}) 
  \quad \text{s.t.} \quad \frac{c(\mathbf{x})}{\varepsilon_c} -1 \leq 0.
\end{equation}

The normalization aids the convergence of the optimizer, which is now discussed.\\

In contrast to the commonly used Method of Moving Asymptotes (MMA), this work employs the Globally Convergent Method of Moving Asymptotes (GCMMA) to ensure convergence to a Karush-Kuhn-Tucker (KKT) point for any initial guess~\cite{Svanberg1987,Svanberg2002}.\
Although the KKT residual does not converge to zero, numerical experimentation showed more stable convergence compared to MMA.\
The MMA parameters are those suggested by Svanberg~\cite{Svanberg1987}: $c=1000$, $a_0$ and all $d_i$ are set to 1 and $a_i=0$.\
The parameters $s_{\mathrm{init}} = 0.2$,  $s_{\mathrm{decr}} = 0.65$ and $s_{\mathrm{incr}} = 1.07$ control the asymptotes.\
For stability, an outer move limit $\mu = 0.05$ is used.\
At the start of a new continuation run, the optimized design variables of the previous run are used as initial guess.\
Additionally, the GCMMA internal state variables, namely the lower and upper asymptotes and the design variable history, are also passed to GCMMA.\\

The convergence criterion is five consecutive iterations where the constraints are satisfied and either the KKT condition or the function condition holds.\
The KKT condition considers a result converged when the norm of the KKT residual is either lower than $1\mathrm{e}-6$, lower than $1\mathrm{e}-4$ for fifteen consecutive iterations or its change is lower than $1\mathrm{e}-7$.\
The function objective considers a result converged if, for ten consecutive iterations, the relative change in the compliance and volume is lower than $1\mathrm{e}-4$.\

\section{Hole nucleation}\label{app:hole_nucleation}

This appendix briefly discusses details regarding the hole nucleation scheme.\\

Given the filter radius $R$ and the $\eta$ threshold values of the projection, the length scale $r_{\min}$ is computed with the formula from Qian and Sigmund~\cite{qian2013topological}.\
Whenever a design is punctured at element $i$, all active elements whose center lies within a radius of $r_{\min}$ have their density value set to zero.\
The puncture is only accepted if this does not change the design's topology, which is computed via the floodfill algorithm described below.\\

\begin{figure}[h]
\centering
\begin{subfigure}[t]{0.5\linewidth}
  \centering
  \includegraphics[width=\textwidth]{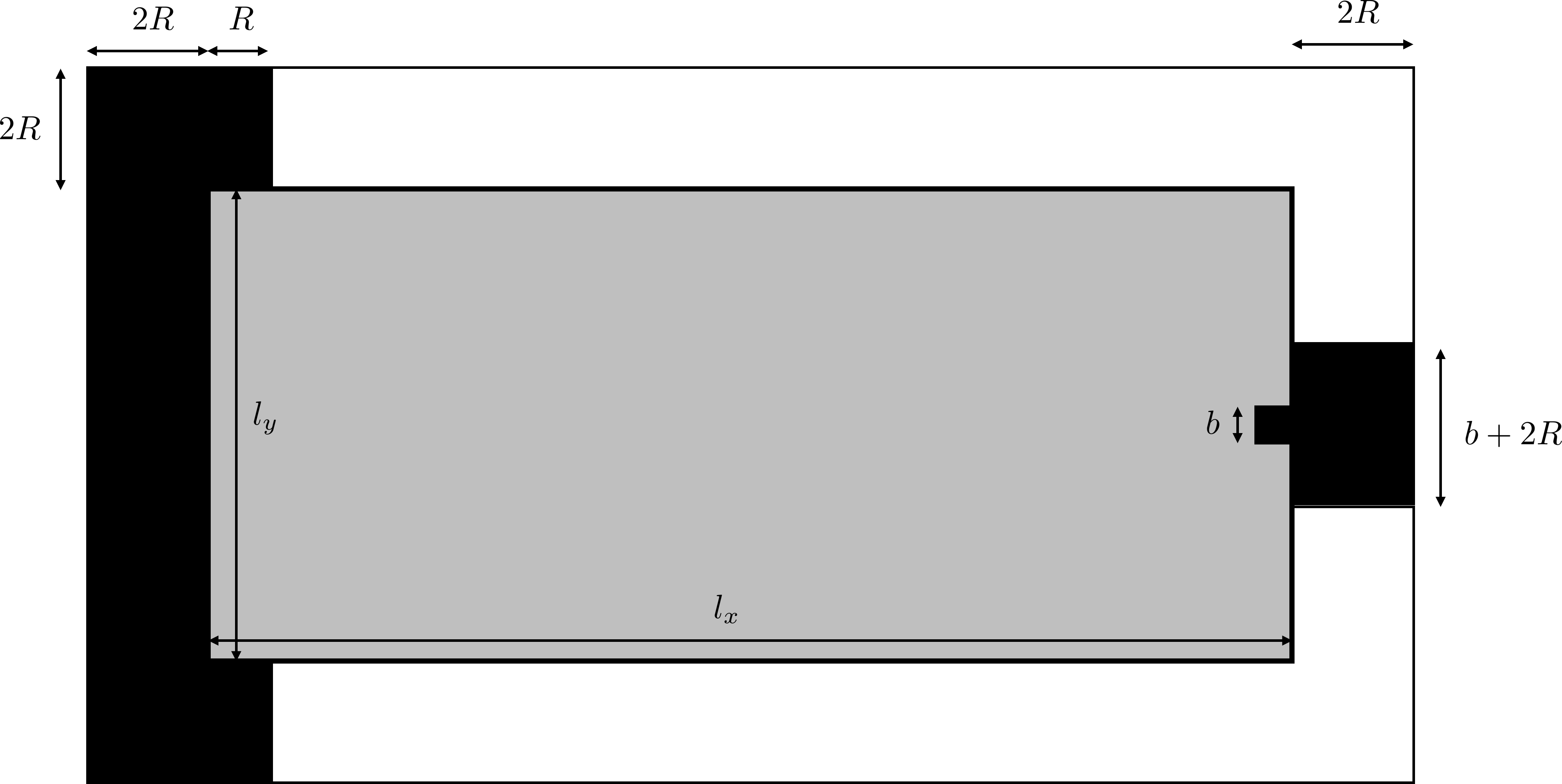}
\caption{Design domain extended with solid at loaded boundaries (+$R$) and void at free boundaries.}
\end{subfigure}%
\begin{subfigure}[t]{0.5\linewidth}
  \centering
  \includegraphics[width=\textwidth]{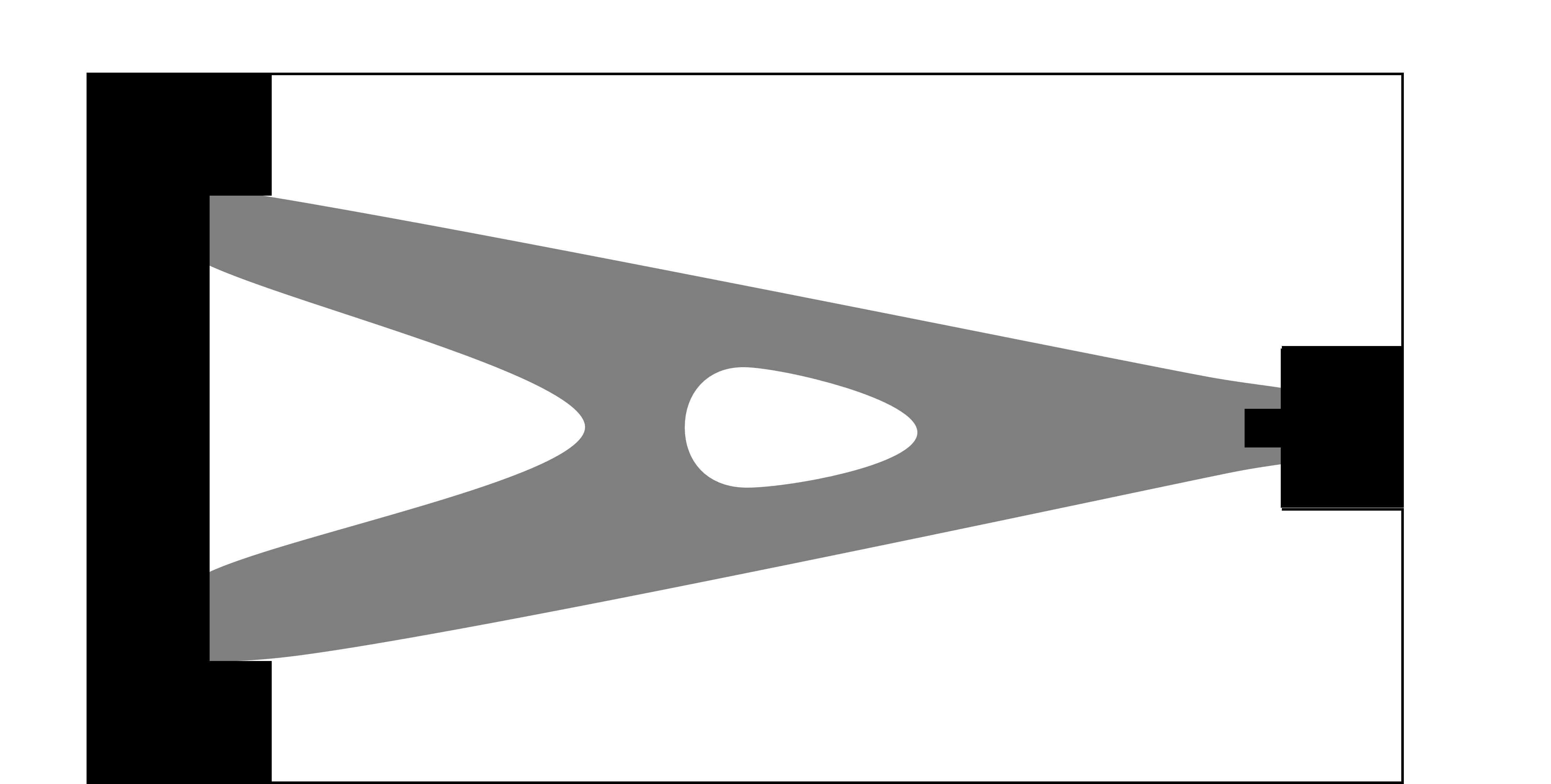}
\caption{Example design with two internal holes.}
\end{subfigure}\\
\begin{subfigure}[t]{0.5\linewidth}
  \centering
  \includegraphics[width=\textwidth]{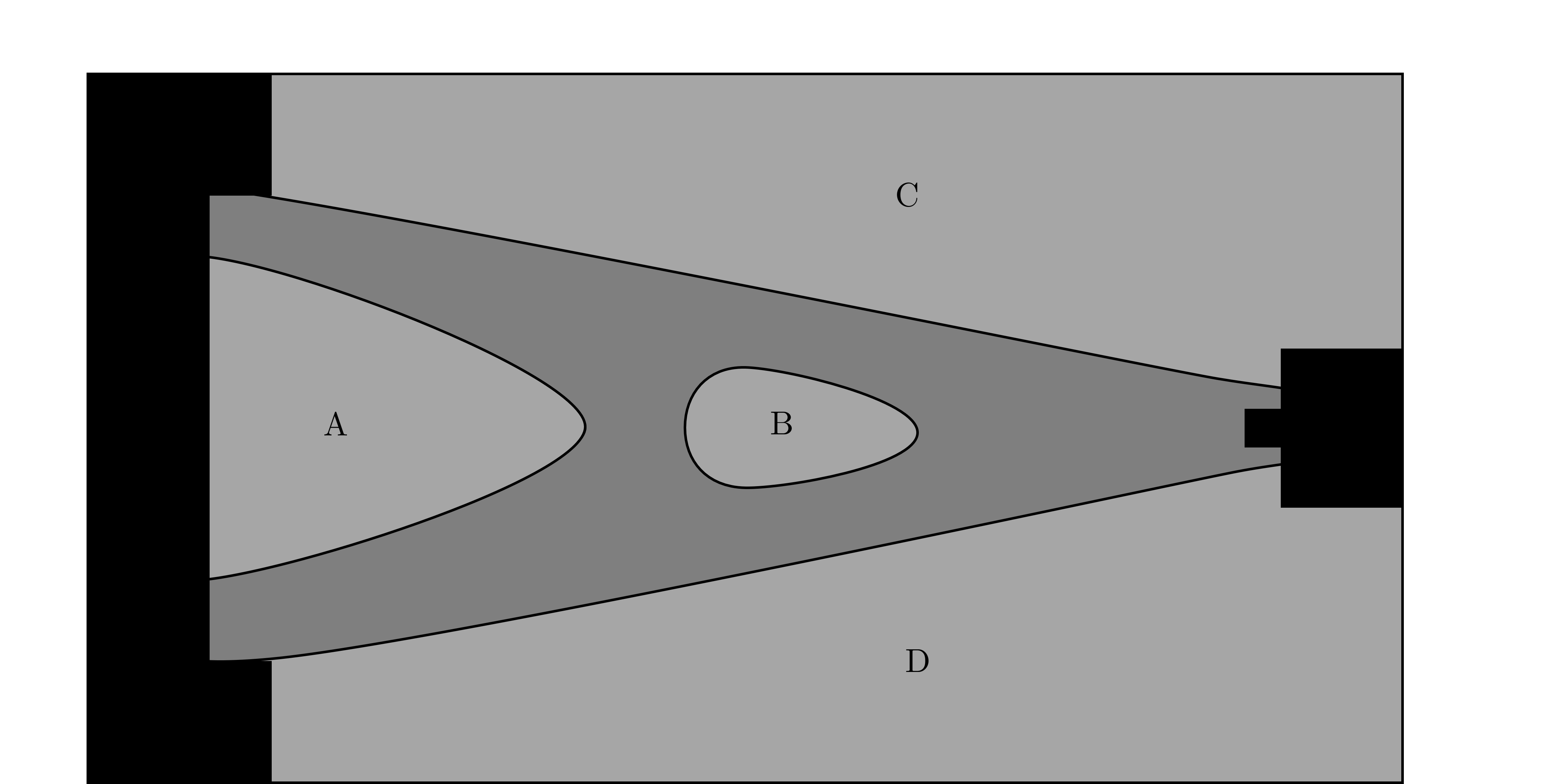}
\caption{Result of hole filling and counting algorithm. The two outer holes ($C$ and $D$) are not counted towards the genus.}
\end{subfigure}%
\caption{Determination of the genus of a design.}
\label{fig:imfill}
\end{figure}

\Cref{fig:imfill} illustrates the way the genus (or topology) of a design is counted.\
The design domain is first extended.\
A solid layer of width $2R$, with $R$ the filter radius, is added to the boundaries on which Robin boundary conditions are applied during the filtering.\
Specifically, these are the boundaries within a distance $R$ of the loaded and fixed boundary.\
All other regions are extended with a void layer of equal width.\
Then, the design itself is projected to solid and void: the blueprint design is mapped to $0$ or $1$ with a cutoff at $0.5$.\
Next, the genus $g$ is counted iteratively by starting at $g=0$ and using Matlab's \texttt{imfill} routine to fill a hole and increment $g$.\
If there are no more holes left to fill, the hole counting routine ends and $g-2$ is returned.\
The $-2$ offset accounts for the two external holes but does not influence the mapping procedure.\\

The hole candidates are distributed throughout the design domain as illustrated in \Cref{fig:hole_grid}.\
This is done by iterating through the elements, starting from the center left and visiting first those to the right and before considering elements above and below the symmetry line.\
This is done to ensure the grid is symmetric.\
For each element, a hole is placed and only accepted as the center of a hole candidate when the resulting hole does not overlap with any of the accepted holes.\
 This condition is verified by checking if the distance between the hole centers is larger than $2r_{\min}+4h$, where the four extra element widths act as a necessary buffer against the erosion operation.\
That is, if designs are punctured with holes closer together, then they are separated by a grey region in the eroded design, which the optimizer can easily remove.\
This changes the topology and leads to premature branch termination.\

\section{Continuation} \label{app:continuation}

This appendix briefly discusses details regarding the hole nucleation scheme.\\

As stated, the continuation scheme chooses $\varepsilon_c$ and $\varepsilon_v$ values from an $n_v \,\times \, n_c$ grid on the $[0,1] \, \times \, [1, c_{\max}]$ region of the $(v(\mathbf{x}), c(\mathbf{x}))$ space.\
In practice, all examples use $n_c=n_v$.\
Successive runs always increase the constraint bounds $\varepsilon_c$ and $\varepsilon_v$ to ensure that every optimization starts from a feasible initial guess.\
This improves numerical stability, but convergence issues are still observed close to $(1, 1)$ due to an extreme flatness of the optimization landscape.\
The flatness is caused by the top and bottom right corners of the cantilever, which contribute very little to the compliance.\
To circumvent this issue, $c_1$ is set to $1.05$ and the first branch is started from a fully solid design $\mathbf{x}_0^\star=\mathbf{1}$.\
The result of the first optimizaton run is denoted $(v_1^V, c_1^V)$ and used to rediscretize the volume by replacing $v_{n_v}=1$ with $v_1^V$.\
Thus, all subsequent branches stop at $v_1^V < 1$.\\

A subcontinuation terminates when (a) the final $\varepsilon$ value is reached, (b) a topological change is detected with the floodfill algorithm detailed in \ref{app:hole_nucleation} or (c) the optimized design is similar to an already obtained design.\
The last case, termed branch collision, is detected via the method described below.\\

Two designs $A$ and $B$, with respective design vectors $\mathbf{x}_A$ and $\mathbf{x}_B$, are considered equal if they meet three conditions.\
First, their compliance must be within a tolerance $\mathrm{tol}_c=0.01$, i.e., $\abs{c(\mathbf{x}_A)-c(\mathbf{x}_B)} \leq \mathrm{tol}_c$.\
Second, the volume must also be within a tolerance $\mathrm{tol}_v=0.001$, i.e., $\abs{v(\mathbf{x}_A)-v(\mathbf{x}_B)} \leq \mathrm{tol}_v$.\
Third, the design vectors themselves should be close as well: $\norm{\mathbf{x}_{b,A} - \mathbf{x}_{b,B}}_1 / n_{\mathrm{el}} \leq \mathrm{tol}_x$, where the subscript $b$ refers to the blueprint realization, $n_{\mathrm{el}}$ is the number of elements and $\mathrm{tol}_x=0.01$.\\

Note that the right value for the tolerances $\mathrm{tol}_c, \mathrm{tol}_v$ and $ \mathrm{tol}_x$ depends on the local flatness of the optimization landscape.\
Because of this complexity, the paper proposes further research into deflation~\cite{Papadopoulos2021} and higher-order continuation~\cite{Cesarano2025,Gangl2026} with Newton-based optimizers, since these can track the curvature of the optimum and adaptively decide closeness of solutions.\ 

\section{Parameters and omitted branches} \label{app:values}

The parameters for the numerical example in \Cref{sec:numerical_examples_1} are listed in \Cref{tab:values}.\
In this table, $n_x$ and $n_y$ are the number of elements in $x$ and $y$, $R$ is the filter radius, $\max_{\mathrm{branch}}$ is the maximum number of iterations of the mapping methodology and $n_c$ and $n_v$ determine the number of $\varepsilon$ values.\
The value for $\mathrm{tol}_{\tau}$ determines the maximal value of the normalized strain energy $\tau$, discussed in \Cref{sec:hole_nucleation}, for an element to be considered for hole nucleation.\
Finally, $n_H$ is the maximum number of holes that can be punctured to spawn a new branch.\
This should be larger than $1$ when symmetry is involved, to allow for the generation of symmetric designs as an initial guess for a new branch.\\

\begin{table}[h!]
\caption{Parameters for numerical examples in \Cref{sec:numerical_examples_1}. The unit for the filter radius $R$ is element widths.}
\label{tab:values}
\begin{tabular}{l c c c c c c c}
\toprule
\textbf{Example} & $n_x$ & $n_y$ & $R$ &  $\max_{\mathrm{branch}}$ & $n_c,n_v$ & $\mathrm{tol}_{\tau}$ & $n_H$ \\
\midrule
Cantilever beam, \Cref{sec:cantilever_beam_example} & $200$ & $100$ & $25$ &  $20$ & $100$ & $0.01$ & $2$ \\[0.5ex]
High-compliance region, \Cref{sec:highC} &  $200$ & $100$ & $25$ &  $4$ & $200$ & $0.01$ & $2$  \\[0.5ex]
Symmetric $2\,\times\,1$ cantilever, \Cref{sec:symmetry_study_2x1} & $200$ & $100$ & $15$ &  $40$ & $200$ & $0.01$ & $4$  \\[0.5ex]
Symmetric $4\,\times\,1$ cantilever, \Cref{sec:symmetry_study_4x1} & $400$ & $100$ & $15$ &  $100$ & $200$  & $0.02$ & $4$ \\[0.5ex]
\bottomrule
\end{tabular}
\end{table}

The flatness study in \Cref{sec:flatness} inherits the relevant parameters (mesh, filtering, ...) from the symmetric $2\,\times\,1$ cantilever study of \Cref{sec:symmetry_study_2x1}.\
The bifurcation study in \Cref{sec:bifurcations} inherits the same parameters but refines the grid used to choose the $\varepsilon_c$ and $\varepsilon_v$ values.\
Namely, the volume minimization runs A-D and E-H divide $[c_0, c_{\max}]$ into one hundred values for $\varepsilon_c$,with $c_{\max}=6$ and $c_0$ the compliance of the initial designs A and E, respectively.\
 The compliance minimization runs I-J divides the volume range $[v_I,v_J] \approx [0.28, 0.36]$ into one hundred values for $\varepsilon_v$,with $v_I$ and $v_J$ the volumes of the designs I and J, respectively.\\
 
The mapping methodology always runs for at most $\max_{\mathrm{branch}}$ iterations, but premature convergence can occur if no eligible hole combinations remain.\
As illustrated in \Cref{fig:workflow}, each iteration spawns a branch and uses continuation to extend it.\
This does not always result in interesting or relevant results, however.\
For completeness, the following discusses each iteration of the mapping procedure for each of the numerical examples in \Cref{sec:numerical_examples_1}.\
More information can be provided upon reasonable request.\\

\Cref{sec:cantilever_beam_example} runs the mapping procedure for $\max_{\mathrm{branch}}=20$ iterations but converged after $7$ iterations due to a lack of eligible hole candidates.\
As it is mainly meant as an illustration, only three branches are shown.\
These correspond to iterations 1, 2 and 5, because iterations 3 and 4 converged to the two-bar truss (found in iteration 2) and iterations 6 an 7 were considered superfluous for the illustrative purpose of the section.\\

\Cref{sec:highC} studies the high compliance region with $\max_{\mathrm{branch}}=4$, for both the symmetry and asymmetric cantilever boundary conditions.\
The visualized branches, namely the 1- and 2-bar trusses, correspond to iterations 1 and 2 for both cases.\
Iterations $3$ and $4$  were not visualized to retain the focus on the crossing of the first two branches.\\

\Cref{sec:symmetry_study_2x1} investigated the $2\,\times\,1$ cantilever with symmetric boundary conditions for $\max_{\mathrm{branch}}=40$, but convergence was reached after iteration $22$ due to a lack of eligible hole candidates.\
The seven branches correspond to iterations 1, 2, 3, 5, 6, 8 and 16.\
The remaining fifteen iterations were not shown because they either (a) did not result in new solutions due to immediate branch collision, (b) converged to previously found branches but did not trigger the equality checks of \ref{app:continuation}, or (c) converged to a mirrored variant of the asymmetric six-bar truss of iteration $5$.\\

\Cref{sec:symmetry_study_4x1} investigated the $4\,\times\,1$ cantilever with symmetric boundary conditions for the maximum of $\max_{\mathrm{branch}}=100$ iterations.\
The fifteen visualized branches correspond to iterations 1, 2, 3, 6, 7, 29, 64, 67, 73, 75, 80, 82, 84, 87 and 95.\
The remaining 85 iterations were not visualized for the same general reasons as for the $2 \, \times\,1$ cantilever.\
Namely, 57 iterations immediately collided with previous branches and hence produced no new solutions, whereas the other 28 iterations yielded (mirrored versions of) previous branches without triggering the equality checks of \ref{app:continuation}.\

\clearpage
\bibliography{mybibfile}

\end{document}